\documentclass[12pt]{article}
\usepackage{latexsym}
\usepackage{amsfonts}

\input amssym.def
\input amssym.tex

\usepackage{euscript}

\begin{document}


\title{Refinements of the Morse stratification of the normsquare of the moment map}

\author{Frances Kirwan\\Mathematical Institute, Oxford University,\\ Oxford OX1 3LB, UK}

 


\newtheorem{prop}{Proposition}[section]
\newtheorem{lem}[prop]{Lemma}
\newtheorem{cor}[prop]{Corollary}
\newtheorem{thm}[prop]{Theorem}
\newtheorem{guess}{Conjecture}
\newtheorem{dff}[prop]{Definition}
\newenvironment{df}{\begin{dff}\rm}{\end{dff}}
\newtheorem{REM}[prop]{Remark}
\newenvironment{rem}{\begin{REM}\rm}{\end{REM}}
\newtheorem{examplit}[prop]{Example}
\newenvironment{example}{\begin{examplit}\rm}{\end{examplit}}

\newcommand{\nn}{\tilde{n}}
\newcommand{\DD}{\tilde{d}}
\newcommand{\mnd}{{\cal M}(n,d)}
\newcommand{\Mnd}{{\cal M}(n,d)}
\newcommand{\dai}{\frac{\partial}{\partial a_{i}}}
\newcommand{\daij}{\frac{\partial^{2}}{\partial a_{i} \partial a_{j}}}
\newcommand{\liek}{{\frak k}}
\newcommand{\lieks}{{\frak k}^*}
\newcommand{\HS}{H^{*}}
\newcommand{\ar}{a_{r}}
\newcommand{\bjr}{b_{r}^{j}}
\newcommand{\fr}{f_{r}}
\newcommand{\Res}{\mbox{\rm Res}}
\newcommand{\ch}{\mbox{\rm ch}}
\newcommand{\mtc}{M_{\T}(c)}
\newcommand{\n}{n}
\newcommand{\D}{d}
\newcommand{\E}{\hat{E}}
\newcommand{\V}{\hat{V}}
\newcommand{\G}{{\cal G}}
\newcommand{\GG}{\overline{\cal G}}
\newcommand{\pa}{\overline{\partial}}
\newcommand{\g}{\bar{g}}
\newcommand{\Cmu}{{\cal C}_{\tau}}
\newcommand{\HG}{H^{*}_{\cal G}}
\newcommand{\Css}{{\cal C}^{ss}}
\newcommand{\C}{{\cal C}}
\newcommand{\liet}{{\frak t}}
\newcommand{\liets}{{\frak t}^*}
\newcommand{\ug}{\tilde{\gamma}}
\newcommand{\ub}{\tilde{\beta}}
\newcommand{\um}{\tilde{\mu}}
\newcommand{\ct}{\tilde{{\cal C}}}
\newcommand{\stab}{{\rm Stab}}

\newcommand{\aff}{{\Bbb A }}
\newcommand{\RR}{{\Bbb R }}
\newcommand{\CC}{{\Bbb C }}
\newcommand{\ZZ}{{\Bbb Z }}
\newcommand{\PP}{ {\Bbb P } }
\newcommand{\QQ}{{\Bbb Q }}
\newcommand{\UU}{{\Bbb U }}

\newcommand{\cala}{{\mbox{$\mathcal A$}}}
\newcommand{\calb}{{\mbox{$\mathcal B$}}}
\newcommand{\calc}{{\mbox{$\mathcal C$}}}
\newcommand{\cald}{{\mbox{$\mathcal D$}}}
\newcommand{\cale}{{\mbox{$\mathcal E$}}}
\newcommand{\calf}{{\mbox{$\mathcal F$}}}
\newcommand{\calg}{{\mbox{$\mathcal G$}}}
\newcommand{\calh}{{\mbox{$\mathcal H$}}}
\newcommand{\cali}{{\mbox{$\mathcal I$}}}
\newcommand{\calj}{{\mbox{$\mathcal J$}}}
\newcommand{\calk}{{\mbox{$\mathcal K$}}}
\newcommand{\call}{{\mbox{$\mathcal L$}}}
\newcommand{\calm}{{\mbox{$\mathcal M$}}}
\newcommand{\caln}{{\mbox{$\mathcal N$}}}
\newcommand{\calo}{{\mbox{$\mathcal O$}}}
\newcommand{\calp}{{\mbox{$\mathcal P$}}}
\newcommand{\calq}{{\mbox{$\mathcal Q$}}}
\newcommand{\calr}{{\mbox{$\mathcal R$}}}
\newcommand{\cals}{{\mbox{$\mathcal S$}}}
\newcommand{\calt}{{\mbox{$\mathcal T$}}}
\newcommand{\calu}{{\mbox{$\mathcal U$}}}
\newcommand{\calv}{{\mbox{$\mathcal V$}}}
\newcommand{\calw}{{\mbox{$\mathcal W$}}}
\newcommand{\calx}{{\mbox{$\mathcal X$}}}
\newcommand{\caly}{{\mbox{$\mathcal Y$}}}
\newcommand{\calz}{{\mbox{$\mathcal Z$}}}


\maketitle

\renewcommand{\theequation}{\thesection.\arabic{equation}}
\newcommand{\renorm}{{ \setcounter{equation}{0} }}


Let $X$ be any nonsingular complex projective variety with a linear action of a 
complex reductive group $G$, and let $X^{ss}$ and $X^s$ be the sets of 
semistable and stable points of $X$ in the sense of Mumford's geometric
invariant theory \cite{MFK}.  We can choose a maximal compact subgroup
$K$ of $G$ and an inner product on the Lie algebra $\liek$ of $K$ which is
invariant under the adjoint action. Then $X$ has a $G$-equivariantly perfect
stratification $\{S_\beta : \beta \in \calb \}$ 
by locally closed nonsingular $G$-invariant subvarieties
with $X^{ss}$ as an open stratum, which can be obtained as the
Morse stratification for the normsquare of a moment map $\mu:X \to \lieks$
for the action of $K$ on $X$ \cite{K2}. In this note 
the Morse stratification $\{S_\beta: \beta \in \calb \}$ 
is refined to obtain stratifications of $X$ by locally closed
nonsingular $G$-invariant subvarieties with $X^s$ as an open stratum. 
The strata can be defined inductively in terms of the sets of stable
points of closed nonsingular subvarieties of $X$ acted on by reductive
subgroups of $G$, and their projectivised normal bundles.

These refinements of the Morse stratification are not in general equivariantly
perfect; that is, the associated equivariant Morse inequalities are not
necessarily equalities. However when $G$ is abelian we can modify
the moment map (or equivalently modify the linearisation of the
action) by the addition of any constant, since the
adjoint action is trivial. Perturbation of the moment map by
adding a small constant then provides an equivariantly perfect
refinement of the stratification $\{S_\beta: \beta \in \calb \}$,
and a generic perturbation gives us a refined stratification whose
strata can be described inductively in terms of the sets of stable
points of linear actions of reductive subgroups of $G$ for which
stability and semistability coincide. This is useful even when
$G$ is not abelian, since important questions about the cohomology of the
 Marsden-Weinstein reduction $\mu^{-1}(0)/K$ (or equivalently
the geometric invariant theoretic quotient $X/\!/G$) can be reduced
to questions about the quotient of $X$ by a maximal torus of $G$.

The same constructions work when $X$ is a compact K\"{a}hler manifold
with a Hamiltonian action of $K$. Even when $X$ is symplectic but
not K\"{a}hler, refinements of the Morse stratification for $|\! | \mu |\! |^2$ 
can be obtained by choosing a suitable almost complex structure and
Riemannian metric.

The moduli spaces $\Mnd$ of holomorphic bundles of rank $n$ and
degree $d$ over a Riemann surface $\Sigma$ of genus $g \geq 2$
can be constructed as quotients of
infinite dimensional spaces of connections, in a way
which is analogous to the construction of quotients in geometric
invariant theory; the r\^{o}le of the moment map is played
by curvature and the r\^{o}le of 
the normsquare of the moment map is played by the Yang-Mills
functional. In \cite{INI} refinements of the
Morse stratification of the Yang-Mills functional are studied
using the ideas of this paper. The
motivation for this study was the search for a complete set of relations
among the standard generators for the cohomology of these moduli
spaces $\Mnd$ when $n$ and $d$ are coprime and $n>2$ \cite{EK}.

The layout of this paper is as follows. $\S$1 reviews some background
material and $\S$2 uses the partial desingularisation construction
of  \cite{K4} to define a stratification $\{ \Sigma_\gamma : \gamma \in \Gamma \}$
of $X^{ss}$ with $X^s$ as an open stratum. $\S$3 gives an inductive description
of the strata $\Sigma_\gamma$ in $X^{ss}$ in terms of the stable and 
semistable points of linear actions of reductive subgroups and subquotients of $G$ on 
nonsingular subvarieties of $X$ and their projectivised normal bundles.
$\S$4 refines the stratification $\{ \Sigma_\gamma : \gamma \in \Gamma \}$
to obtain strata which are described inductively purely in terms of the
stable points (not the semistable points) of the linear actions appearing
in $\S$3, and in $\S$5 this stratification is used to refine the Morse
stratification $\{ S_\beta: \beta \in \calb\}$ of $X$. In $\S$6 an alternative refinement
of this stratification is obtained when $G$ is abelian by perturbing the
moment map; in this case the inner product on the Lie algebra of $K$ can also be perturbed
to give a refined stratification. Finally $\S$7 discusses applications to the
study of the cohomology ring of the
 Marsden-Weinstein reduction $\mu^{-1}(0)/K$ (or equivalently
the geometric invariant theoretic quotient $X/\!/G$), and in particular the relationship
between this cohomology ring and the corresponding cohomology ring when $G$
is replaced by a maximal torus.

\section{Background}
\renorm

In this section we shall review briefly the material we need from
\cite{K2} (see also \cite{G2,MFK,N2}.

Let $X$ be a connected nonsingular projective variety embedded in 
projective space $\PP_n$ and let $G$ be a complex
reductive group acting linearly on $X$ via a homomorphism
$\rho :G \to GL(n+1;\CC)$. Then $G$ is the complexification of
a maximal compact subgroup $K$, and by rechoosing the
coordinates on $\PP_n$ if necessary, we can assume that
$K$ acts unitarily on $\PP_n$ via $\rho: K \to U(n+1)$.

The geometric invariant theoretic quotient $X/\!/G$ 
is the projective variety whose homogeneous
coordinate ring is the $G$-invariant part of the homogeneous
coordinate ring of $X$.
A point $x$ of $X$ is called semistable if there exists an
invariant homogeneous polynomial $f$ which does not vanish on $x$,
and $x$ is called stable if in addition the orbit $Gx$ is a closed
subset of the set of stable points $X^{ss}$ and has dimension $\dim G$. 
There is a rational map from $X$ to $X/\!/G$ which restricts to a
$G$-invariant surjective morphism
$$\phi:  X^{ss} \to X/\!/G$$
from the open subset $X^{ss}$ of $X$, and
every fibre of $\phi$ which meets the set $X^s$ consisting of the points
of $X$ which are stable for the action of $G$ is a single $G$-orbit
in $X^{ss}$. The image of $X^s$ is an open subset of $X/\!/G$
which can thus be identified via $\phi$ with the quotient $X^s/G$. 
We shall assume that $X^s \neq \emptyset$ (but see Remark \ref{rem2.1}).

$X$ has a K\"{a}hler structure given by the restriction of the
Fubini-Study metric on $\PP_n$, and the K\"{a}hler form
$\omega$ is a $K$-invariant symplectic form in $X$. There is
a moment map $\mu :X \to \lieks$, where $\liek$ is the Lie
algebra of $K$, defined by
$$\mu(x) = \rho^*((2 \pi i |\!|x^*|\!|^2)^{-1} x^*  \bar{x}^{*t})$$
for $x \in X \subseteq \PP_n$ represented by $x^*\in \CC^{n+1} \backslash \{0\}$,
when the Lie algebra of $U(n+1)$ and its dual are both
identified with the space of skew-Hermitian $(n+1)\times (n+1)$
matrices in the standard way.
Then $\mu^{-1}(0) $ is a subset of $X^{ss}$, and the inclusion
induces a homeomorphism from the Marsden-Weinstein
reduction (or symplectic quotient) $\mu^{-1}(0)/K$ to the
geometric invariant theoretic quotient $X/\!/G$. In fact
$X^{ss}$ is the set of points $x$ in $X$ such that $\mu^{-1}(0)$
meets the closure of the $G$-orbit of $x$, and $X^s$ is the set of 
points $x$ in $X$ such that $\mu^{-1}(0)$ meets the $G$-orbit of $x$ 
in a point which is regular for $\mu$. In the good case when $X^{ss} = X^s$
then $X/\!/G = X^{ss}/G$ and its rational cohomology is isomorphic to
the equivariant cohomology of $X^{ss}$.

If we fix an invariant inner product on the Lie algebra
$\liek$ then we can consider the function $|\!|\mu|\!|^2$
as a Morse function on $X$. It is not in general a Morse function in
the classical sense, nor even a Morse-Bott function, since
the connected components of its set of critical points may
not be submanifolds of $X$, but nonetheless it induces
a Morse stratification $\{ S_{\beta} : \beta \in {\cal B}\}$
of $X$ such that the stratum to which $x \in X$ belongs is determined by
the limit set of its path of steepest descent for $|\!|\mu|\!|^2$ with
respect to the Fubini-Study metric. This stratification can also
be defined purely algebraically and has the following properties \cite{K2}.

\begin{prop} \label{prop1.1}  i) Each stratum $S_{\beta}$ is a $G$-invariant
locally closed nonsingular subvariety of $X$.

\noindent ii) The unique open stratum $S_0$ is the set $X^{ss}$ of
semistable points of $X$.

\noindent iii) The stratification is equivariantly perfect over the rationals, so that
$$\dim H^i_G(X) = \dim H^i_G(X^{ss}) + \sum_{\beta \neq 0}
\dim H^{i-\lambda(\beta)}_G(S_{\beta})$$
where $\lambda(\beta)$ is the real codimension of $S_{\beta}$
in $X$.

\noindent iv) If $\beta \neq 0$ then there is a proper nonsingular 
subvariety $Z_{\beta}$ of $X$ acted on by a reductive subgroup
$\mbox{\stab }\beta$ of $G$ such that
$$H^*_G(S_{\beta}) \cong H^*_{\mbox{\stab }(\beta)}(Z_{\beta}^{ss})$$
where $Z_{\beta}^{ss}$ is the set of semistable points of $Z_{\beta}$
with respect to an appropriate linearisation of the action of $\mbox{\stab }(\beta)$.

\end{prop}

\begin{rem}
i) All cohomology in this paper has rational coefficients.

\noindent ii) We shall assume that the invariant inner product chosen on $\liek$
is rational, and we shall use it to identify $\lieks$ with $\liek$ throughout.

\noindent iii) Note that $G$-equivariant cohomology is the same as $K$-equivariant cohomology, since
$G$ retracts onto its maximal compact subgroup $K$.
\end{rem}

If we choose a positive Weyl chamber $\liet_+$ in the Lie
algebra $\liet$ of a maximal torus $T$ of $K$, then we can
identify the indexing set ${\cal B}$ with a finite subset of $\liet_+$
(or equivalently with the set of adjoint orbits of the points
in this finite subset of $\liet_+$, since each adjoint orbit meets
$\liet_+$ in exactly one point).
More precisely, let $\alpha_0,...,\alpha_n$ be the weights
of the representation of $T$ on $\CC^{n+1}$ and use the 
restriction to $\liet$ of the fixed invariant inner product on
$\liek$ to identify $\liets$ with $\liet$. Then $\beta \in \liet_+$
belongs to ${\cal B}$ if and only if $\beta$ is the closest point
to 0 of the convex hull in $\liet$ of some nonempty subset of
$\{ \alpha_0,...,\alpha_n \}$. Moreover $\mbox{\stab }(\beta)$ is the stabiliser
of $\beta$ under the adjoint action of $G$, and $Z_{\beta}$ is
the intersection of $X$ with the linear subspace 
$$\{ [x_0:...:x_n] \in \PP_n  : x_i =0 \mbox{ if } \alpha_i . \beta  \neq |\!| \beta |\!|^2 \}$$
of $\PP_n$. Equivalently $Z_{\beta}$ is the union of those connected
components of the fixed point set of the subtorus $T_{\beta}$ of $T$
generated by $\beta$ on which the constant value taken by the
real-valued function $x \mapsto \mu(x).\beta$ is $|\! | \beta |\! |^2$. 
If we use the fixed inner product to identify $\liek$ with its dual,
then the image $\mu(Z_{\beta})$ of $Z_{\beta}$ is contained in
the Lie algebra of the maximal compact subgroup $\mbox{\stab }_K(\beta) = K \cap \mbox{\stab }(\beta)$
of $\mbox{\stab }(\beta)$. Since moment maps are only unique up to the addition of
a central constant, we can take $\mu - \beta$ as our moment map for the
action of $\mbox{\stab }_K(\beta)$ on $Z_{\beta}$. This corresponds to a modification of
the linearisation of the action of $\mbox{\stab }(\beta)$ on $Z_{\beta}$ whose restriction to the
complexification $T_{\beta}^c$ of $T_{\beta}$ is trivial, and we define
$Z_{\beta}^{ss}$ to be the set of semistable points of $Z_{\beta}$ with respect
to this modified linear action. Equivalently, $Z_{\beta}^{ss}$ is the stratum
labelled by 0 (the minimum stratum)  for the Morse stratification of the function $|\!|\mu - \beta|\!|^2$ on
$Z_{\beta}$. Then
\begin{equation} \label{sb}
S_{\beta} = G \bar{Y_{\beta}} \backslash \bigcup_{|\! | \gamma |\! | > |\! | \beta |\! |} G \bar{Y_{\gamma}}
=G Y_{\beta}^{ss} \end{equation}
where
$$\bar{Y}_{\beta} = \{ x \in X : x_i = 0 \mbox{ if } \alpha_i . \beta < |\!| \beta |\!|^2 \}$$
and 
$$Y_{\beta} = \{ x \in \bar{Y}_{\beta} : x_i \neq \mbox{0 for some $i$ such that }\alpha_i . \beta = |\!| \beta |\!|^2 \},$$
while 
\begin{equation} \label{yb} Y_{\beta}^{ss} = p_{\beta}^{-1}(Z_{\beta}^{ss}) \end{equation}
where $p_{\beta}:Y_{\beta} \to Z_{\beta}$ is the obvious projection
given by 
\begin{equation} \label{pb} p_{\beta}(x) = \lim_{t \to \infty} \exp (-it\beta) x. \end{equation}

\begin{rem} \label{zbss} In \cite{K2} $Z_\beta^{ss}$ is defined as above to be the set of semistable points
of $Z_\beta$ with respect to the  action of $\stab (\beta)$ with the linearisation modified in a
way which corresponds to replacing the moment map $\mu$ by (a positive integer
multiple of) $\mu - \beta$. However
we can, if we wish,  regard $\mu - \beta$ as a moment map for the action of
$\stab _K (\beta)/T_{\beta}$ on $Z_{\beta}$, and then $Z_{\beta}^{ss}$ is the set of semistable points for the
induced linear action of its complexification $\stab (\beta)/T_{\beta}^c$. The advantage of this description
is that the analogous definition of $Z_\beta^s$ is a useful one, whereas $Z_\beta$ has no points which
are stable with respect to the action of $\stab(\beta)$ when $\beta \neq 0$. Therefore we shall define
$Z_\beta^s$ to be the set of stable points for the action of $\stab (\beta)/T^c_\beta$ on $Z_\beta$,
linearised so that the corresponding moment map is a positive integer multiple of $\mu - \beta$.
Equivalently
$Z_\beta^{ss}$ and $Z_\beta^s$ are the sets of semistable points of $Z_\beta$
under the action of the subgroup of $\stab(\beta)$ whose Lie algebra is the complexification
of the orthogonal complement to $\beta$ in the Lie algebra of $\stab_K(\beta)$.

\end{rem}

In fact if $B$ is the Borel
subgroup of $G$ associated to the choice of positive Weyl chamber $\liet_+$ and if
$P_{\beta}$ is the parabolic subgroup $B \mbox{\stab }(\beta)$, then $Y_{\beta}$ and
$Y_{\beta}^{ss}$ are $P_{\beta}$-invariant and we have
\begin{equation} \label{SGY} S_{\beta} \cong G \times_{P_{\beta}} Y_{\beta}^{ss}. \end{equation}
Moreover $Y_{\beta}$ is a nonsingular subvariety of $X$ and $p_{\beta}:Y_{\beta}
\to Z_{\beta}$ is a locally trivial fibration whose fibre is isomorphic to
$\CC^{m_{\beta}}$ for some $m_{\beta} \geq 0$.

Now an element $g$ of $G$ lies in the parabolic subgroup $P_{\beta}$ if and only if
$\exp(-it\beta) g \exp(it\beta)$ tends to a limit in $G$ as $t \to \infty$, and this limit defines
a surjection $q_{\beta}: P_{\beta} \to \stab (\beta)$.
Since the surjections $p_{\beta}:Y_{\beta}^{ss} \to Z_{\beta}^{ss}$ and 
$q_{\beta}:P_{\beta} \to \mbox{\stab } (\beta)$ 
are retractions
satisfying 
\begin{equation} \label{pqp} p_{\beta}(gx) = q_{\beta}(g) p_{\beta}(x) \end{equation}
for all $g \in P_{\beta}$ and $x \in Y_{\beta}^{ss}$, this gives us
the isomorphism
$$H^*_G(S_{\beta}) \cong H^*_{\mbox{\stab }(\beta)}(Z_{\beta}^{ss})$$
of Proposition \ref{prop1.1}(iv). Moreover since $G=KB$ and $B \subseteq P_{\beta}$
we have $G\bar{Y_{\beta}} = K \bar{Y_{\beta}}$, which is compact, and hence
\begin{equation} \label{closure} \bar{S_{\beta}} \subseteq G\bar{Y_{\beta}}
\subseteq S_{\beta} \cup \bigcup_{|\!| \gamma |\!| > |\!| \beta |\!|} S_{\gamma}. \end{equation}

Note that $x = [x_0,...,x_n] \in X$ is semistable (respectively stable) for the action of the complex torus
$T^c$ if and only if 0 belongs to the convex hull in $\liet$ of the set of weights
$\{ \alpha_j:x_j \neq 0 \}$ (respectively to the interior of this convex hull), and that
$x$ is semistable (respectively stable) for the action of $G$ if and only if every
element $gx$ of its $G$-orbit is semistable (respectively stable) for the action of $T^c$.
In particular this tells us that
$$S_0 = X^{ss}.$$
It also tells us that if $\beta \in \calb$ 
 and if $y \in \bar{Y}_\beta \backslash Y^{ss}_\beta,$
then there is a subset $S$ of $\{\alpha_j:\alpha_j . \beta \geq |\!|\beta |\!|^2 \}$ such that
\begin{equation} \label{ch} y \in \stab(\beta) Y_{\beta'} \end{equation}
where $\beta'$ is the closest point to 0 of the convex hull of $S$
and $\beta' \neq \beta$.

\begin{rem} \label{rhonot} Of course the definitions of 
the subvarieties $Z_{\beta}$, $Z_{\beta}^{ss}$ and
so on in this section depend upon the action of $G$ on $X \subseteq \PP_n$ and its
linearisation via the representation $\rho: G \to GL(n+1;\CC)$. When it is necessary to
make these explicit in the notation we shall write $Z_{\beta}(X, \rho)$, $Z_{\beta}^{ss}(X, \rho)$ and
so on, but we shall omit $X$ if $X=\PP_n$, so that for example 
$$Z_{\beta}(X,\rho) = X \cap Z_{\beta}(\rho).$$
\end{rem}

\begin{rem} In the special case when every semistable point is stable,
the quotient variety $X/\!/G$ is topologically the ordinary quotient
$X^{ss}/G$ where $G$ acts with only finite stabilisers on $X^{ss}$,
which means that its Betti numbers $\dim H^i(X/\!/G)$ are the
same as the equivariant Betti numbers $\dim H^i_G(X^{ss})$ of
$X^{ss}$. These can be calculated inductively using Proposition
\ref{prop1.1} together with the fact that the equivariant cohomology of a nonsingular complex
projective variety is isomorphic as a vector space to the tensor product of
its ordinary cohomology and the equivariant cohomology of
a point.
\end{rem}

\begin{rem}
If $X$ is any compact symplectic manifold with a Hamiltonian action of a
compact Lie group $K$ and a moment map $\mu:X \to \lieks$, then the
Morse stratification of $|\!| \mu |\!|^2$ has essentially the properties described
here, with $G$ replaced by $K$ and $P_\beta$ replaced by $K \cap \stab(\beta)$
\cite{K2}.
\end{rem}

\section{Stratifying the set of semistable points}

\renorm

Suppose now that $X$ has some stable points but also has semistable
points which are not stable. 
In \cite{K4,K6} it is described how one can blow $X$ up along a sequence of 
nonsingular  $G$-invariant subvarieties to obtain a $G$-invariant morphism 
$\tilde{X} \to X$ where $\tilde{X}$ is a complex projective variety acted on 
linearly by $G$ such that $\tilde{X}^{ss} = \tilde{X}^s$. The induced birational 
morphism $\tilde{X}/\!/G \to X/\!/G$ of the geometric invariant theoretic quotients 
is then a partial desingularisation of $X/\!/G$ in the sense that $\tilde{X}/\!/G$ 
has only orbifold singularities (it is locally isomorphic to the quotient of a 
nonsingular variety by a finite group action) whereas the singularities of $X/\!/G$ 
are in general much worse. In this section we shall review the construction
of the partial desingularisation $\tilde{X}/\!/G$ and use it firstly to stratify $X^{ss}$
and subsequently (in $\S$5) to refine the stratification $\{ S_{\beta}:\beta \in {\cal B} \}$ of $X$
described in $\S$1.

The set $\tilde{X}^{ss}$ can be obtained from $X^{ss}$ as follows. There exist 
semistable points of $X$ which are not stable if and only if there exists a 
non-trivial connected reductive subgroup of $G$ fixing a semistable point. Let 
$r>0$ be the maximal dimension of a reductive subgroup of $G$ 
fixing a point of $X^{ss}$ and ${\cal R}(r)$ be a set of representatives of conjugacy 
classes of all connected reductive subgroups $R$ of 
dimension $r$ in $G$ such that 
$$ Z^{ss}_{R} = \{ x \in X^{ss} :  \mbox{$R$ fixes $x$}\} $$
is non-empty. Then
$$
\bigcup_{R \in {\cal R}(r)} GZ^{ss}_{R}
$$
is a disjoint union of nonsingular closed subvarieties of $X^{ss}$. The action of 
$G$ on $X^{ss}$ lifts to an action on the blow-up $X_{(1)}$ of 
$X^{ss}$ along $\bigcup_{R \in {\cal R}(r)} GZ_R^{ss}$ which can be linearised so that the complement 
of $X_{(1)}^{ss}$ in $X_{(1)}$ is the proper transform of the 
subset $\phi^{-1}(\phi(GZ_R^{ss}))$ of $X^{ss}$ where $\phi:X^{ss} \to X/\!/G$ is the quotient 
map (see \cite{K4} 7.17). Moreover no point of $X_{(1)}^{ss}$ is fixed by a 
reductive subgroup of $G$ of dimension at least $r$, and a point in $ X_{(1)}^{ss} $ 
is fixed by a reductive subgroup $R$ of 
dimension less than $r$ in $G$ if and only if it belongs to the proper transform of the 
subvariety $Z_R^{ss}$ of $X^{ss}$.

We can now apply the same procedure to $X_{(1)}^{ss}$ to obtain $X_{(2)}^{ss}$ such 
that no reductive subgroup of $G$ of dimension at least 
$r-1$ fixes a point of $X_{(2)}^{ss}$. If we repeat this process enough times, we obtain 
$X_{(0)}^{ss} = X^{ss},X_{(1)}^{ss},X_{(2)}^{ss},\ldots,X_{(r)}^{ss}$ such that
no reductive subgroup of $G$ of positive dimension fixes a point of $X_{(r)}^{ss}$,
and we set $\tilde{X}^{ss} = X_{(r)}^{ss}$. Equivalently we can construct a sequence  
$$X_{(R_0)}^{ss} = X^{ss}, X_{(R_1)}^{ss},\ldots,X_{(R_\tau)}^{ss} = \tilde{X}^{ss}$$ 
where $R_1,\ldots,R_\tau$ are connected reductive subgroups of $G$ with 
$$r= \dim R_1 \geq \dim R_2 \geq \cdots \dim R_\tau \geq 1,$$
 and if $1 \leq l \leq \tau$ then
 $X_{(R_l)}$ is the blow up of $X_{(R_{l-1})}^{ss}$ 
along its closed nonsingular subvariety $GZ_{R_l}^{ss}\cong G \times_{N_l} Z_{R_l}^{ss}$,
where $N_l$ is the normaliser of $R_l$ in $G$. Similarly 
$\tilde{X}/\!/G = \tilde{X}^{ss}/G$ can be obtained from $X/\!/G$ by blowing 
up along the proper transforms of the images $Z_R /\!/N$ 
in $X/\!/G$ of the subvarieties $GZ_R^{ss}$ of $X^{ss}$ in decreasing order of $\dim R$.

If $1 \leq l \leq \tau$ then there is a $G$-equivariant stratification 
$$\{ \cals_{\beta,l}: \beta \in {\cal B}_l \}$$ 
of $X_{(R_l)}$
by nonsingular 
$G$-invariant locally closed subvarieties such that one of the strata, indexed by $0  \in {\cal B}_l$, 
coincides with the open 
subset $X_{R_l}^{ss}$ of $X_{R_l}$.  This stratification is constructed exactly as the stratification
$\{ S_{\beta} : \beta \in {\cal B} \}$ of $X$ was constructed in the last section; note that $X_{(R_l)}$ is
in general only quasi-projective rather than projective, but it is shown in \cite{K4} that the
construction of the stratification
still works for $X_{(R_l)}$ and the properties given in Proposition \ref{prop1.1} still hold. 

There is a partial ordering on ${\cal B}_l$ with $0 $ as its 
minimal element such that if 
$\beta \in {\cal B}_l$ then the closure in $X_{(R_l)}$ of the stratum $\cals_{\beta,l}$ satisfies
$$
\overline{\cals_{\beta,l}} \subseteq \bigcup_{\gamma \in {\cal B}_l, \gamma \geq \beta} \cals_{\gamma,l}.
$$
If $\beta \in {\cal B}_l$ and $\beta \neq 0 $ then the stratum $\cals_{\beta,l}$ retracts $G$-equivariantly 
onto its (tranverse) 
intersection with the exceptional divisor $E_l$ for the blow-up $X_{(R_l)} \to X_{(R_{l-1})}^{ss}$. This 
exceptional divisor is isomorphic to the projective bundle $\PP({\cal N}_l)$ over $G\hat{Z}_{R_l}^{ss}$,
where $\hat{Z}_{R_l}^{ss}$ is the proper transform of $Z_{R_l}^{ss}$ in $X_{(R_{l-1})}^{ss}$ and
${\cal N}_l$ is 
the normal bundle to $G \hat{Z}_{R_l}^{ss}$ in $X_{R_{l-1}}^{ss}$. The stratification 
$\{\cals_{\beta,l}: \beta \in {\cal B}_l\}$ is determined 
by the action of $R_l$ on the fibres of ${\cal N}_l$ over $Z_{R_l}^{ss}$ (see \cite{K4} \S 7).

The composition 
$$
\tilde{X}^{ss} = X^{ss}_{(R_\tau)} \to X^{ss}_{(R_{\tau-1})} \to \cdots X^{ss}_{(R_1)} \to X^{ss}
$$
is an isomorphism over the set  $X^s$ of stable points of
 $X$, and the complement of $X^s$ in $\tilde{X}^{ss}$  is 
just the union of the proper transforms in $\tilde{X}^{ss}$ 
of the exceptional divisors $E_1,\ldots,E_k$ for the blow-ups $X_{R_l} \to X^{ss}_{(R_{l-1})}$ for $l=1,\ldots,\tau$.

We can now stratify $X^{ss}$ as follows. We take as the highest stratum 
the nonsingular closed subvariety $GZ_{R_1}^{ss}$ whose 
complement in $X^{ss}$ can be naturally identified with 
the complement $X_{(R_1)}\backslash E_1$ of the exceptional divisor
$E_1$ in $X_{(R_1)}$.
Recall that $GZ^{ss}_{R_1} \cong G \times_{N_1} Z^{ss}_{R_1}$ where $N_1$ is the
normaliser of $R_1$ in $G$, and $Z^{ss}_{R_1}$ is equal to the set of semistable
points for the action of $N_1$, or equivalently for the induced action of $N_1/R_1$,
on $Z_{R_1}$, which is a union of connected components of the fixed point set of
$R_1$ in $X$ (see \cite{K4} \S 5). Since $R_1$ has maximal dimension among those
reductive subgroups of $G$ with fixed points in $X^{ss}$, we have
$$Z_{R_1}^{ss} = Z_{R_1}^s$$
where $Z_{R_l}^s$ is the set of stable points for the action of $N_l/R_l$
on $Z_{R_l}$ for $1 \leq l \leq \tau$.

Next we take as strata the nonsingular locally closed subvarieties
$$
\{ \cals_{\beta,1}\backslash E_1: \beta \in {\cal B}_1, \beta \neq 0\}
$$
of $X_{(R_1)}\backslash E_1 = X^{ss}\backslash GZ_{R_1}^{ss}$, whose complement 
in $X_{(R_1)}\backslash E_1$ is just 
$X_{(R_1)}^{ss}\backslash E_1=X_{(R_1)}^{ss}\backslash E_1^{ss}$ where $E_1^{ss} = X_{(R_1)}^{ss} \cap E_1$,
and then we take the intersection of 
$X_{(R_1)}^{ss}\backslash E_1$ with $GZ_{R_2}^{ss}$. This intersection is $G Z_{R_2}^s$ where
$Z_{R_2}^s$ is the set of stable points for the action of $N_2/R_2$ on $Z_{R_2}$, and its
complement in $X_{(R_1)}^{ss}\backslash E_1$ can be naturally identified with the 
complement in $X_{(R_2)}$ of the union of $E_2$ and the 
proper transform $\hat{E}_1$ of $E_1$.

Our next strata are the nonsingular locally closed subvarieties
$$
\{ \cals_{\beta,2} \backslash (E_2 \cup \hat{E}_1): \beta \in {\cal B}_2, \beta \neq 0\}
$$
of $X_{(R_2)}\backslash(E_2 \cup \hat{E}_1)$, whose complement in $X_{(R_2)}\backslash(E_2 \cup \hat{E}_1)$ is 
$X_{(R_2)}^{ss}\backslash(E_2 \cup \hat{E}_1)$, and the stratum after these 
is $GZ_{R_3}^s$. Repeating this process gives us strata which are all nonsingular locally 
closed $G$-invariant subvarieties of $X^{ss}$ indexed by the disjoint union
$$
 \{R_1\} \, \cup \, \{ R_1\} \! \times \! ({\cal B}_1 \backslash \{0\}) \, \cup \, \cdots  
\, \cup \,  \{R_\tau\} \, \cup \, \{R_\tau\} \! \times \! ({\cal B}_\tau \backslash \{0\}), $$
and the complement in $X^{ss}$ of the union of these strata is just the open subset $X^s$.
We take $X^s$ as our final stratum indexed by 0, so that the indexing set for our
stratification of $X^{ss}$ is the disjoint union
\begin{equation}
\label{new*}
{\Gamma} = 
 \{R_1\} \, \cup \, \{ R_1\} \times ({\cal B}_1 \backslash \{0\}) \, \cup \, \cdots  
\, \cup \,  \{R_\tau\} \, \cup \, \{R_\tau\}\times({\cal B}_\tau \backslash \{0\}) \, \cup \, \{ 0 \}.
\end{equation}
 Moreover the given partial 
orderings on ${\cal B}_1,\ldots,{\cal B}_\tau$ together with the ordering in the expression (\ref{new*}) above
for $\Gamma$
induce a partial ordering on $\Gamma$, with $R_1$ as the maximal element
and 0 as the minimal element, such that the closure in $X^{ss}$ of the 
stratum $\Sigma_\gamma$ indexed by $\gamma \in \Gamma$ satisfies
\begin{equation} \label{stratclos}
\overline{\Sigma_\gamma} \subseteq \bigcup_{\delta \in \Gamma, \delta \geq \gamma } \Sigma_{\delta}.
\end{equation}
Thus this process gives us a stratification 
\begin{equation} \label{strat1} \{ \Sigma_{\gamma}: \gamma \in \Gamma \} \end{equation}
of $X^{ss}$
such that 
the stratum indexed by the minimal element 0 of $\Gamma$ coincides with the open subset $X^s$ of
$X^{ss}$. 

\begin{rem} \label{rem2.1}
We have been assuming that $X^s \neq \emptyset$, but this procedure gives us a stratification of
$X^{ss}$ even when $X^s$ is empty. The only difference when $X^s$ is empty is that the procedure
terminates at some stage $l$ when
$$X^{ss}_{(R_{l-1})} = G \hat{Z}^{ss}_{R_l} \cong  G \times_{N_l} \hat{Z}^{ss}_{R_l}$$
and gives us a stratification indexed by
$${\Gamma} = 
 \{R_1\} \, \cup \, \{ R_1\} \times ({\cal B}_1 \backslash \{0\}) \, \cup \, \cdots  
\, \cup \,  \{R_{l-1}\} \, \cup \, \{R_{l-1}\}\times({\cal B}_{l-1} \backslash \{0\}) \, \cup \, \{ R_l \}
$$
such that the stratum indexed by the minimal element $R_l$ of $\Gamma$ is the open subset
$GZ_{R_l}^s$ of $X^{ss}$. Note also that $Z_{R_l}^s$ is nonempty, since otherwise
$$Z_{R_l}^{ss} = N_l Z_R^{ss}$$
for some $R$ containing $R_l$ with $\dim R > \dim R_l$, and then
$$GZ_R^{ss} = GZ_R^{ss} = X^{ss},$$
so the procedure would have terminated at an earlier stage.
\end{rem}

\section{Inductive description of the strata $\Sigma_{\gamma}$ in $X^{ss}$}

\renorm

The last section described  a stratification 
$ \{ \Sigma_{\gamma}: \gamma \in \Gamma \} $
of $X^{ss}$
such that 
the stratum indexed by the minimal element 0 of $\Gamma$ coincides with the open subset $X^s$ of
$X^{ss}$. 
In this section we shall study the strata $\Sigma_{\gamma}$ in more detail.

Note that the strata $\Sigma_{\gamma}$ with $\gamma \neq 0$ fall into two classes. Either
$\gamma = R_l$ for some $l \in \{1,...,\tau\}$, in which case the stratum $\Sigma_\gamma$ is
$$GZ^s_{R_l},$$
or else $\gamma =(R_l,\beta)$ where $\beta \in {\cal B}_l \backslash \{ 0 \}$ for some $l \in \{1,...,\tau\}$ and
the stratum $\Sigma_{\gamma}$ is
$$\cals_{\beta,l} \backslash (E_l \cup \hat{E}_{l-1} \cup ... \cup \hat{E}_1).$$
In the latter case we know from (\ref{SGY}) that
\begin{equation} \label{2.2} \cals_{\beta,l} = G Y^{ss}_{\beta,l}
 \cong G \times_{P_{\beta,l}} Y^{ss}_{\beta,l} \end{equation}
where $Y^{ss}_{\beta,l}$ fibres over $Z^{ss}_{\beta,l}$ via
$p_{\beta}:Y^{ss}_{\beta,l} \to Z^{ss}_{\beta,l}$ with fibre $\CC^{m_{\beta,l}}$
for some $m_{\beta,l} >0$, and
\begin{equation} \label{2.3} \cals_{\beta,l}\cap E_l = G (Y^{ss}_{\beta,l}\cap E_l) \cong G 
\times_{P_{\beta,l}} (Y^{ss}_{\beta,l} \cap E_l) \end{equation}
where $Y^{ss}_{\beta,l}\cap E_l$ fibres over $Z^{ss}_{\beta,l}$ with fibre
$\CC^{m_{\beta,l}-1}$ (see \cite{K4} Lemmas 7.6 and 7.11). Thus
\begin{equation} \label{2.4} \cals_{\beta,l} \backslash E_l  \cong 
G \times_{P_{\beta,l}} (Y^{ss}_{\beta,l} \backslash E_l)
\end{equation}
where $Y^{ss}_{\beta,l} \backslash E_l$ fibres over $Z^{ss}_{\beta,l}$ with fibre
$\CC^{m_{\beta,l}-1} \times (\CC \backslash \{ 0 \})$. Let 
$$\pi_l: E_l \cong \PP({\cal N}_l) \to G \hat{Z}^{ss}_{R_l}$$
denote the projection. Lemma 7.9 of \cite{K4} tells us that if $x \in \hat{Z}^{ss}_{R_l}$ then
the intersection of $\cals_{\beta,l}$ with the fibre $\pi_l^{-1}(x) = \PP({\cal N}_{l,x})$
of $\pi_l$ at $x$ is the union of  those strata indexed by points in the adjoint orbit
$Ad(G)\beta$ in the stratification of $\PP({\cal N}_{l,x})$ induced by the representation
$\rho_l$ of $R_l$ on the normal ${\cal N}_{l,x}$ to $G\hat{Z}^{ss}_{R_l}$ at $x$.
Note that we can assume that $R_l \cap K$ is a maximal compact subgroup of
$R_l$ and that $R_l\cap T$ is a maximal torus for $R_l \cap K$, and then
$Ad(G)\beta$ meets a positive Weyl chamber for $R_l$ in $\mbox{Lie}(R_l \cap T)$ in
a finite number of points 
$$\beta = \beta_1 = Ad(w_1)\beta, \; \beta_2=Ad(w_2)
\beta,\ldots ,\beta_{r_{\beta,l}}=Ad(w_{r_{\beta,l}})\beta$$ 
where $w_1=1,w_2,\ldots,  w_{r_{\beta,l}} \in
G$ represent elements of the Weyl group of $G$.

Now if $y \in Z^{ss}_{\beta,l}$ and $\pi_l(y) = gx$ where $g \in G$ and $x \in \hat{Z}_{R_l}^{ss}$,
then $x$ is fixed by $Ad(g^{-1})\beta$ and so $Ad(g^{-1})\beta$ lies in the Lie algebra of
$R_l$. Since $R_l \cap T$ is a maximal compact torus of $R_l$, there exists $r \in R_l$ such that 
$Ad(rg^{-1})\beta \in Lie(R_l \cap T)$. Then $Ad(rg^{-1})\beta = Ad(w_j) \beta$ for some 
$j \in \{1,...,r_{\beta,l} \}$, and hence $w_j^{-1} r g^{-1} \in \stab (\beta)$, so $g \in \stab (\beta)
 w_j^{-1} R_l$. Conversely if $g=hw_j^{-1} r$ where $h \in \stab (\beta)$ and $r \in R_l$, then $y \in Z^{ss}_{\beta,
l}$ if and only if $h^{-1} y $ lies in $Z^{ss}_{\beta,l}$ where
$\pi_l (h^{-1} y) = w_j^{-1} r x \in w_j^{-1} \hat{Z}_{R_l}^{ss}$. Thus 
$$Z_{\beta,l}^{ss} = \bigcup_{1 \leq j \leq r_{\beta,l}} \stab (\beta) \left( 
Z_{\beta,l}^{ss} \cap w_j^{-1} \pi_l^{-1}(\hat{Z}_{R_l}^{ss})  \right)$$
$$ = \bigcup_{1 \leq j \leq r_{\beta,l}} \stab (\beta) w_j^{-1} \left( 
Z_{Ad(w_j)\beta,l}^{ss} \cap  \pi_l^{-1}(\hat{Z}_{R_l}^{ss})  \right).$$
Also if $\stab (\beta) ( Z_{\beta,l}^{ss} \cap w_j^{-1} \pi_l^{-1}(\hat{Z}_{R_l}^{ss}) )$ meets
$\stab (\beta) ( Z_{\beta,l}^{ss} \cap w_i^{-1} \pi_l^{-1}(\hat{Z}_{R_l}^{ss}) )$ then
$\stab (\beta) w_j^{-1} \hat{Z}_{R_l}^{ss}$ meets $w_i^{-1} \hat{Z}_{R_l}^{ss}$, and since
$G\hat{Z}_{R_l}^{ss} \cong G \times_{N_l} \hat{Z}_{R_l}^{ss}$ this means that there is some
$h \in \stab (\beta)$ and $n \in N_l$ such that $w_i h = n w_j$, so that 
$\beta_i = Ad(w_i) \beta \in Ad(N_l) \beta_j$.  Conversely if $\beta_i \in Ad(N_l) \beta_j$
then $w_i h = n w_j$ for some $h \in \stab (\beta)$ and $n \in N_l$, and so
$$\stab (\beta) ( Z_{\beta,l}^{ss} \cap w_j^{-1} \pi_l^{-1}(\hat{Z}_{R_l}^{ss}) )
= \stab (\beta) ( Z_{\beta,l}^{ss} \cap h^{-1}w_i^{-1} n \pi_l^{-1}(\hat{Z}_{R_l}^{ss}) )$$
$$= \stab (\beta) ( Z_{\beta,l}^{ss} \cap w_i^{-1} \pi_l^{-1}(\hat{Z}_{R_l}^{ss}) ).$$
Thus $Z_{\beta,l}^{ss}$ is a disjoint union
$$Z_{\beta,l}^{ss} = \bigsqcup_{1 \leq j \leq s_{\beta,l}} \stab (\beta)  \left( 
Z_{\beta,l}^{ss} \cap w_j^{-1} \pi_l^{-1}(\hat{Z}_{R_l}^{ss})  \right)$$
where $Ad(w_1)\beta = \beta,\ldots, Ad(w_{s_{\beta,l}})\beta$ form a set
of representatives for the $Ad(N_l)$ orbits in $Ad(G)\beta$, and $Y^{ss}_{\beta,l}$ and
$\cals_{\beta,l}$ can be expressed similarly as disjoint unions. In fact, since by 
(\ref{pqp}) the fibration
$$p_{\beta}: Y^{ss}_{\beta,l} \to Z^{ss}_{\beta,l}$$
satisfies $p_\beta(gy) = q_\beta(g) p_\beta(y)$ for all $g \in P_\beta$ and
$y \in Y^{ss}_{\beta,l}$ where $q_\beta : P_\beta \to \stab(\beta)$ is the
projection, we have
$$Y_{\beta,l}^{ss} = \bigsqcup_{1 \leq j \leq s_{\beta,l}} P_\beta p_\beta^{-1}  \left( 
Z_{\beta,l}^{ss} \cap w_j^{-1} \pi_l^{-1}(\hat{Z}_{R_l}^{ss})  \right)$$
and
$$\cals_{\beta,l} = \bigsqcup_{1 \leq j \leq s_{\beta,l}} G p_\beta^{-1}  \left( 
Z_{\beta,l}^{ss} \cap w_j^{-1} \pi_l^{-1}(\hat{Z}_{R_l}^{ss})  \right).$$

This means that we could, if we wished, replace the indexing set ${\cal B}_l \backslash \{0  \}$, whose elements
correspond to the $G$-adjoint orbits $Ad(G)\beta$ of elements of the indexing set for the
stratification of $\PP({\cal N}_{l,x})$ induced by the representation $\rho_l$, by the set of their
$N_l$-adjoint orbits $Ad(N_l)\beta$. Then we would still have (\ref{2.2}) -- (\ref{2.4}), but now
\begin{equation} \label{2.5}
Z_{\beta,l}^{ss} = \stab (\beta) ( Z_{\beta,l}^{ss} \cap \pi_l^{-1}(\hat{Z}_{R_l}^{ss})  )
\cong \stab (\beta) \times_{N_l \cap \stab (\beta)} ( Z_{\beta,l}^{ss} \cap \pi_l^{-1}(\hat{Z}_{R_l}^{ss})  )
\end{equation}
and
\begin{equation} \label{new2.6}
Y_{\beta,l}^{ss} \cong P_\beta \times_{Q_\beta} p_\beta^{-1}  \left( 
Z_{\beta,l}^{ss} \cap w_j^{-1} \pi_l^{-1}(\hat{Z}_{R_l}^{ss})  \right) \end{equation}
where
$$Q_\beta = q_\beta^{-1} (N_l \cap \stab(\beta))$$
is a subgroup of $P_\beta$, and hence
\begin{equation} \label{new2.7}
\cals_{\beta,l} \cong G \times_{P_\beta} Y_{\beta,l}^{ss} \cong
G \times_{Q_\beta} p_\beta^{-1}  \left( 
Z_{\beta,l}^{ss} \cap w_j^{-1} \pi_l^{-1}(\hat{Z}_{R_l}^{ss})  \right). \end{equation}
Furthermore $\pi_l$ now restricts to a fibration
\begin{equation} \pi_l:Z_{\beta,l}^{ss} \cap \pi_l^{-1}(\hat{Z}_{R_l}^{ss})  )
\to \hat{Z}_{R_l}^{ss} \end{equation}
whose fibre at $x \in \hat{Z}_{R_l}^{ss}$ is $Z_{\beta}^{ss}(\rho_l)$ defined
as at Remark \ref{rhonot}, where $\rho_l$ is the representation of $R_l$ on the normal ${\cal N}_{l,x}$
to $G \hat{Z}_{R_l}^{ss}$ at $x$.

If $1 \leq j \leq l-1$ then the proper transform $\hat{E}_j$ in $X^{ss}_{(R_l)}$ of the exceptional
divisor $E_j$ in $X^{ss}_{(R_j)}$ meets the exceptional divisor $E_l \cong \PP({\cal N}_l)$ transversely,
and their intersection is the restriction 
$$\PP({\cal N}_l |_{\hat{E}_j \cap G\hat{Z}_{R_l}})$$
of the projective bundle $\PP({\cal N}_l)$ over $G\hat{Z}^{ss}_{R_l}$ to the intersection in
$X^{ss}_{(R_{l-1})}$ of $G\hat{Z}_{R_l}^{ss}$ with the proper transform of $E_j$ in 
$X^{ss}_{(R_{l-1})}$ (which by abuse of notation we shall also denote
by $\hat{E}_j$). Moreover the complement in  $G\hat{Z}_{R_l}^{ss}$ of its intersection with
the exceptional divisors $\hat{E}_1,...,\hat{E}_{l-1}$ is just $G {Z}_{R_l}^s$.  Thus
\begin{equation} \label{2.6}
\Sigma_{\gamma} = G Y^{\backslash E}_{\beta,l} \cong 
G \times_{P_{\beta}} Y^{\backslash E}_{\beta,l} \end{equation}
where 
$$Y^{\backslash E}_{\beta,l} = Y^{ss}_{\beta,l} \backslash (E_l \cup \hat{E}_{l-1} \cup ... \cup \hat{E}_1)$$
 fibres over
$\stab (\beta) \times_{N_l \cap \stab (\beta)} ( Z_{\beta,l}^{ss} \cap \pi_l^{-1}(Z_{R_l}^s)  )$
with fibre $\CC^{m_{\beta,l}-1} \times (\CC\backslash \{0\})$, and
$Z_{\beta,l}^{ss} \cap \pi_l^{-1}(Z_{R_l}^s) $ fibres over 
$Z^s_{(R_l)}$ with fibre $Z_{\beta}^{ss}(\rho_l)$. In addition, if we set
$$Y^{\backslash E}_\beta = Y^{\backslash E}_{\beta,l} \cap p_\beta^{-1} \left(
Z^{ss}_{\beta,l} \cap \pi^{-1}_l(\hat{Z}^{ss}_{R_l}) \right)$$
we have from
(\ref{new2.6}) and (\ref{new2.7}) that
\begin{equation} \label{ypqy} Y^{\backslash E}_{\beta,l} \cong P_\beta
\times_{Q_\beta} Y^{\backslash E}_\beta \end{equation}
and hence
\begin{equation} \label{sgqy}
\Sigma_\gamma \cong G \times_{Q_\beta} Y^{\backslash E}_\beta \end{equation}
where $Q_\beta = q_\beta^{-1} (N_l \cap \stab(\beta))$ and
$p_\beta: Y^{\backslash E}_\beta \to Z^{ss}_{\beta,l} \cap \pi^{-1}_l(\hat{Z}^{ss}_{R_l}) $
is a fibration with fibre $\CC^{m_{\beta,l}-1} \times (\CC\backslash \{0\})$.

\begin{rem}
The moduli space
$\mnd$ of semistable holomorphic bundles of rank $n$ and degree $d$ over a fixed Riemann
surface of genus $g \geq 2$ can be constructed as a quotient of an
infinite dimensional affine space of connections ${\cal C}$ by a complexified gauge
group $\G_c$,
in an infinite-dimensional version of the 
construction of quotients in geometric invariant theory, or equivalently as an
infinite dimensional symplectic reduction with curvature as a moment map.
When $n$ and $d$ are coprime, semistability coincides with stability and 
$\mnd$ is the topological quotient of the semistable
subset ${\cal C}^{ss}$ of ${\cal C}$ by the action of $\G_c$.
The r\^{o}le of
the normsquare of the moment map is played by the Yang-Mills
functional, which was studied by Atiyah and Bott in
their fundamental paper \cite{AB}. Atiyah and Bott studied
the stratification of ${\cal C}$ defined using the Harder-Narasimhan
type of a holomorphic bundle, which they expected 
to be the Morse stratification of the Yang-Mills 
functional (this was later shown to be the case \cite{Daskal}).
The methods of this paper can be 
used to provide a stratification $\{\Sigma_\gamma:\gamma \in \Gamma \}$ of
${\cal C}^{ss}$ with ${\cal C}^s$ as the unique open stratum.
This stratification of $\calc^{ss}$ and induced refinements of the
Yang-Mills stratification of $\calc$ are studied in detail in \cite{INI},
where they are related to natural refinements of the notion of the Harder-Narasimhan
type of a holomorphic
bundle.
\end{rem}

\section{A refined stratification of $X^{ss}$}

\renorm

We can now iterate the construction of the stratification (\ref{strat1}) described in 
$\S$2 and use induction on the dimension of $G$ to define a
stratification 
\begin{equation} \label{strat2}
\{ \tilde{\Sigma}_{\ug }:\ug  \in \tilde{\Gamma} \}
\end{equation}
of $X^{ss}$ by $G$-equivariant nonsingular subvarieties
which refines the stratification (\ref{strat1}). When the dimension of $G$ is
zero, so that $X^s = X^{ss} = X$, then
$\tilde{\Gamma} = \Gamma$ and the stratification has one stratum which is $X$ itself. When 
$\dim G >0$ then we shall refine the stratification $\{ \Sigma_{\gamma}: \gamma \in \Gamma \}$
defined at (\ref{strat1}) as follows. 
If $\gamma \in \Gamma \backslash \{0,R_1,...,R_\tau\}$ then $\gamma = (R_l,\beta)$
where $\beta \in {\cal B}_l
\backslash \{0 \}$ for some $l \in \{1,...,\tau\}$, and by (\ref{2.6}) we have
$$
\Sigma_{\gamma} = GY^{\backslash E}_{\beta,l} \cong 
G \times_{P_{\beta}} Y^{\backslash E}_{\beta,l} $$
where $Y^{\backslash E}_{\beta,l}$ fibres over
$\stab (\beta) \times_{N_l \cap \stab (\beta)} ( Z_{\beta,l}^{ss} \cap \pi_l^{-1}(Z_{R_l}^s)  )$
with fibre $\CC^{m_{\beta,l}-1} \times (\CC\setminus \{0\})$, and
$Z_{\beta,l}^{ss} \cap \pi_l^{-1}(Z_{R_l}^s) $ fibres over $Z^s_{(R_l)}$ with fibre $Z_{\beta}^{ss}(\rho_l)$.
We have a linear action of $R_l \cap \stab (\beta) /T^c_\beta$ on $Z_\beta(\rho_l)$ which corresponds
(up to multiplication by a positive integer) to the moment map $\mu - \beta$.
 Therefore by induction on $\dim G$ we can assume that
we have defined a stratification
$\{ \tilde{\Sigma}^{\gamma}_{\ub }: \ub  \in \tilde{\Gamma}_{\gamma}\}$
of $Z_\beta^{ss}(\rho_l)$ by nonsingular $R_l \cap \stab (\beta)$-invariant subvarieties.
In fact, since the stabiliser in $G$ of any $x \in Z_R^s$  has connected
component $R_l$, we can assume that we have a stratification 
of $Z_\beta^{ss}(\rho_l)$ by nonsingular $\stab (x) \cap \stab (\beta)$-invariant subvarieties.
Since $\stab (x) \subseteq N_l$ and since the fibration
$$\pi_l:Z^{ss}_{\beta,l} \cap \pi_l^{-1}(Z^s_{R_l}) \to Z^s_{R_l}$$
is $N_l \cap \stab (\beta)$-equivariant with fibre $Z_\beta^{ss}(\rho_l)$, this gives us
a stratification of $Z^{ss}_{\beta,l} \cap \pi_l^{-1}(Z^s_{R_l})$ by nonsingular
$N_l \cap \stab (\beta)$-invariant subvarieties, and hence a stratification of
$$\stab (\beta) \times_{N_l \cap \stab (\beta)} ( Z_{\beta,l}^{ss} \cap \pi_l^{-1}(Z_{R_l}^s)  )$$
by nonsingular $\stab (\beta)$-invariant subvarieties. We also have a fibration
$$p_\beta:Y^{\backslash E}_{\beta,l} = 
Y^{ss}_{\beta,l} \backslash (E_l \cup \hat{E}_{l-1} \cup ... \cup \hat{E}_1) \to
\stab (\beta) \times_{N_l \cap \stab (\beta)} ( Z_{\beta,l}^{ss} \cap \pi_l^{-1}(Z_{R_l}^s)  )$$
with fibre $\CC^{m_{\beta,l}-1} \times (\CC \backslash \{0\})$, which satisfies $p_\beta(gx) = q_\beta(g) p_\beta(x)$
for all $g \in P_\beta$ and $x \in Y^{\backslash E}_{\beta,l} $
(see (\ref{pqp})). Thus we get an induced stratification of 
$Y^{\backslash E}_{\beta,l} $ by $P_\beta$-invariant
subvarieties, and finally an induced stratification of
$$
\Sigma_{\gamma}  \cong 
G \times_{P_{\beta}} Y^{\backslash E}_{\beta,l} $$
by nonsingular $G$-invariant subvarieties $\tilde{\Sigma}^\gamma_{\ub }$ for
$\ub  \in \tilde{\Gamma}_\gamma$. In particular $\Sigma_{\gamma}$ has
an open stratum
\begin{equation} \label{siggams} \Sigma^s_\gamma = \tilde{\Sigma}^\gamma_0 \end{equation}
corresponding to the open stratum $Z^s_\beta(\rho_l)$ of $Z^{ss}_\beta(\rho_l)$
consisting of stable points for the action of $R_l \cap \stab (\beta) / T_\beta^c$.

In this way we obtain a stratification $\{\tilde{\Sigma}_{\ug}: \ug \in \tilde{\Gamma} \}$
of $X^{ss}$ indexed by 
$$
\tilde{\Gamma} = \{\ug = (\gamma):\gamma \in \{0,R_1,...,R_\tau\} \} 
\cup \{ \ug = (R_l,\beta): 1\leq l \leq \tau \mbox{ and } \beta \in \calb_l \backslash \{ 0 \} \}$$
$$\cup \{ \ug = (R_l,\beta,\gamma_1,
\ldots ,\gamma_t):t\geq 1 \mbox{ and }   1\leq l \leq \tau \mbox{ and } \beta \in \calb_l \backslash \{ 0 \}$$
\begin{equation} \label{2.9} 
\mbox{ and }  (\gamma_1,\ldots,\gamma_t) \in \tilde{\Gamma}_{(R_l,\beta)}\backslash \{ 0 \} \} ,
\end{equation}
where $\tilde{\Gamma}_{(R_l,\beta)}$ is defined inductively as above, 
and the strata $\tilde{\Sigma}_{\ug}$ are given by 
$$\tilde{\Sigma}_{(0)} = X^s$$
and if $1 \leq l \leq \tau$ and $\beta \in \calb_l \backslash \{ 0 \}$ then
$$\tilde{\Sigma}_{(R_l)} = G Z_{R_l}^s \mbox{ and }\tilde{\Sigma}_{(R_l,\beta)} = \Sigma^s_{(R_l,\beta)}  ,$$
while if $\ug = (R_l, \beta, \gamma_1,\ldots,\gamma_t)$  then 
$$\tilde{\Sigma}_{\ug} = \tilde{\Sigma}^{(R_l,\beta)}_{(\gamma_1,\ldots,\gamma_t)}.$$

\section{The refined Morse stratification}

In $\S$2 a stratification $\{ \Sigma_\gamma: \gamma \in \Gamma\}$ of the set
$X^{ss}$ of semistable points of $X$ was defined, and in $\S$4 this stratification
was refined to give  a stratification $\{\tilde{\Sigma}_{\ug}: \ug \in \tilde{\Gamma} \}$
of $X^{ss}$. Via the inductive description of the strata $S_\beta$ of the Morse stratification
of $|\!| \mu |\!|^2$
in terms of the semistable points of nonsingular subvarieties of $X$ given
in Proposition \ref{prop1.1},  we can use these stratifications of $X^{ss}$ to
refine the Morse stratification. 

For simplicity we shall just discuss the
stratification $\{ \tilde{\Sigma}_{\tilde{\gamma}}: \tilde{\gamma} \in \tilde{\Gamma}\}$ of $X^{ss}$ constructed in
$\S$4.
The construction of this stratification can be applied for each $\beta \in {\cal B}$ to the 
action of $\mbox{\stab }(\beta)$ 
 on
the nonsingular projective 
subvariety $Z_{\beta}$ of $X$ which appeared in Proposition \ref{prop1.1}(iv) 
(or more precisely to the action of the quotient $\stab (\beta)/T_{\beta}^c$ of $\stab (\beta)$ by its complex
subtorus $T^c_{\beta}$ which acts trivially on $Z_{\beta}$)
to give a stratification
$$\{ \tilde{\Sigma}^{[\beta]}_{\ug}: \ug \in \tilde{\Gamma}_{[\beta]}\}$$
of $Z_{\beta}^{ss}$ by nonsingular $\stab (\beta)$-invariant subvarieties, with $Z_{\beta}^s$ as the stratum
indexed by $(0)$. 
Since $S_{\beta} = GY^{ss}_{\beta}$ satisfies (\ref{SGY}) and we have a retraction $p_{\beta}:
Y^{ss}_{\beta} \to Z^{ss}_{\beta}$ satisfying (\ref{pqp}), we can stratify $S_{\beta}$ as the disjoint union of
strata
\begin{equation} \label{strat} 
 G p_{\beta}^{-1}(\tilde{\Sigma}^{[\beta]}_{\ug}) \cong 
  G \times_{P_{\beta}} p_{\beta}^{-1}(\tilde{\Sigma}^{[\beta]}_{\ug})
\end{equation}
for $\gamma \in \tilde{\Gamma}_{[\beta]}$. This gives us a stratification
\begin{equation} \label{refstrat} \label{strat3} \{ \tilde{S}_{\ub } : 
\ub  \in \tilde{{\cal B}}  \} \end{equation}
of $X$  indexed by
$$\tilde{{\cal B}} = \tilde{\Gamma} \cup \bigcup_{\beta \in {\cal B}
\backslash \{0\} } \{ \beta \} \times \tilde{\Gamma}_{[\beta]}$$
where $ \tilde{S}_{\ub } = \tilde{\Sigma}_{\ub}$ defined as in $\S$4 if $\ub  \in \tilde{\Gamma}$,
and if $\ub  = (\beta_1,..., \beta_t)$ where $\beta_1 \in {\cal B} \backslash \{0\}$ and
$(\beta_2,...,\beta_t) \in \tilde{\Gamma}_{[\beta]}$ then
$$\tilde{S}_{\ub } = G p_{\beta_1}^{-1}( \tilde{\Sigma}^{[\beta_1]}_{(\beta_2,...,\beta_t)}).$$
This stratification $ \{ \tilde{S}_{\ub } : 
\ub  \in \tilde{{\cal B}}  \}$
refines the original stratification $\{S_{\beta}:\beta \in {\cal B} \}$ and has the following
properties.

\begin{prop} \label{prop2.1}  i) Each stratum $\tilde{S}_{\ub }$ is a $G$-invariant
locally closed nonsingular subvariety of $X$.

\noindent ii) The unique open stratum $\tilde{S}_{(0)}$ is the set $X^{s}$ of
stable points of $X$.

\noindent iii) There is a partial ordering $>$ on $\tilde{{\cal B}}$ such that
if $\ub  \in \tilde{{\cal B}}$ then the closure
$\overline{\tilde{S}_{\ub }}$ in $X$ of the stratum $\tilde{S}_{\ub }$
satisfies
$$\overline{\tilde{S}_{\ub }} \subseteq \bigcup_{\ug \geq \ub }
\tilde{S}_{\ug}.$$

\noindent iv) $\tilde{{\cal B} }$ has a subset $\tilde{\Gamma}$ such that
$\ug < \ub $ for all $\ug \in \tilde{\Gamma}$ and $\ub 
\in \tilde{{\cal B}} \backslash \tilde{\Gamma}$ and
$$\bigcup_{\ug \in \tilde{\Gamma}} \tilde{S}_{\ug} = X^{ss}.$$

\noindent v) If $\ub  \in \tilde{{\cal B}} $ and $\ub  \neq (0)$
then the stratum $\tilde{S}_{\ub }$ can be described inductively
in terms of the sets of stable points of certain nonsingular linear sections $Z$ of $X$
and projectivised normal
bundles of nonsingular subvarieties of $X$, acted on by reductive subgroups of $G$
and their quotients.
\end{prop}

\begin{rem}
If $X$ is a compact K\"{a}hler manifold which has a Hamiltonian action of a compact group
$K$ with moment map $\mu : X \to \lieks$, then the Morse stratification for $|\!| \mu |\!|^2$
can be refined just as in Proposition \ref{prop2.1}. Even when $X$ is symplectic but not
K\"{a}hler we can construct a similar refinement by choosing a suitable almost complex
structure and Riemannian metric on $X$ (cf. \cite{K4,MSj}).
\end{rem}

\begin{example} \label{previo}
Consider the action of $G=SL(2;\CC)$ 
and its maximal compact subgroup $K=SU(2)$ on 
$X=(\PP_1)^n$, with the moment map given by the centre of gravity
in $\RR^3$ when $\PP_1$ is identified suitably with the unit sphere
in $\RR^3$ and $\RR^3$ is identified with the Lie algebra of $SU(2)$.
An element $(x_1,\ldots,x_n)$ of $(\PP_1)^n$ is semistable (respectively 
stable) for the action of $G$ if and
only if at most $n/2$ (respectively strictly fewer
than $n/2$) of the points $x_j$ coincide
anywhere on $\PP_1$. 
The Morse stratification for the normsquare of the moment map on $X$ has strata
$S_0=X^{ss}$ and $S_{2j-n}$ for $n/2 < j \leq n$. If $n/2 < j \leq n$
then the elements of $S_{2j-n}$ correspond to sequences of $n$ points on
$\PP_1$ such that exactly $j$ of these points coincide somewhere on
$\PP_1$, and $S_{2j-n}$ retracts equivariantly onto
the subset of $X$ where $j$ points coincide somewhere on $\PP_1$ and
the remaining $n-j$ points coincide somewhere else on $\PP_1$, which
is a single $G$-orbit with stabilizer $\CC^*$, for $j<n$, and with stabiliser
a Borel subgroup of $G$ when $j=n$  (see \cite{K2} $\S$16.1 for more details).

If $n$ is odd then semistability coincides with stability and the refined
stratification coincides with the Morse stratification of $X$.
Now suppose
that $n$ is even, so that semistability and stability do not coincide. The 
semistable elements of $X$ which are fixed by nontrivial
connected reductive subgroups of $G$ are those represented by
sequences $(x_1,...,x_n)$ of points of $\PP_1$
such that there exist distinct $p$ and $q$ in $\PP_1$  
with exactly half of the points $x_1,...,x_n$ equal to $p$ and the rest equal
to $q$. They form 
$n!/2((n/2)!)^2$ $G$-orbits, and their stabilisers are all 
conjugate to the maximal torus $T_c = \CC^*$ of $G$. 
These stabilisers act with weights 2 and $-2$, each with multiplicity $(n/2) -1$,  on 
the normals to the orbits.  We obtain the 
partial desingularization $\tilde{X}/\!/G$ by blowing up $X/\!/G$ at the points
corresponding to these orbits, or equivalently by blowing up $X^{ss}$
along these orbits, removing from the blowup the unstable points (which
form the proper transform of the set of $(x_1,...,x_n)\in X^{ss}$ such
that exactly half of the points $x_1,...,x_n$ coincide somewhere on $\PP_1$)
and finally quotienting by $G$. 
The refined stratification 
 $ \{ \tilde{S}_{\ub } : 
\ub  \in \tilde{{\cal B}}  \}$
 of $X$ thus has
as its strata the set $\tilde{S}_{(0)} = X^s$ of stable points, the set 
$\tilde{S}_{(T)}$ consisting of points represented by
sequences $(x_1,...,x_n)$ in $\PP_1$
such that there exist distinct $p$ and $q$ in $\PP_1$  
with exactly half of $x_1,...,x_n$ equal to $p$ and the rest equal
to $q$, the set $\tilde{S}_{(T,2)}$ consisting of points represented by
sequences $(x_1,...,x_n)$ in $\PP_1$
such that there exists $p$ in $\PP_1$  
with exactly half of the points $x_1,...,x_n$ equal to $p$ and the rest different
from $p$ and not all equal to each other, and finally the strata $S_{2j-n}$
(for $n/2 < j \leq n$) of the Morse stratification.
\end{example}

\begin{example} \label{example5.2}
A very similar example is given by the action of $G=SL(2;\CC)$ 
on 
$X=\PP_n$ identified with the space of unordered 
sequences of $n$ points in $\PP_1$; that is,
with the projectivized symmetric product $\PP(S^n(\CC^2))$  (see \cite{K2} $\S$16.2 for more details).
The diagonal subgroup $T \cong S^1$ of $K=SU(2)$
acts with weights $n,n-2,n-4,...,2-n,-n$ on $S^n(\CC^2)=\CC^{n+1}$.
An element $[a_0,...,a_n]$ of $\PP_n$ corresponds to
the $n$ roots in $\PP_1$ of the polynomial in one variable $t$, say, with
coefficients $a_0,...,a_n$; it is semistable (respectively 
stable) for the action of $G$ if and
only if at most $n/2$ (respectively strictly fewer
than $n/2$) of these roots coincide
anywhere on $\PP_1$, and 
the Morse stratification for the normsquare of the moment map on $X$ 
is essentially the same as in Example \ref{previo} above.
Again when $n=2m$ is even the 
semistable elements of $X$ which are fixed by nontrivial
connected reductive subgroups of $G$ are those represented by
polynomials 
such that there exist distinct $p$ and $q$ in $\PP_1$  
with exactly half the roots equal to $p$ and the rest equal
to $q$. They form a single $G$-orbit, and the stabiliser $ \CC^*$  
 acts with weights $\pm 4, \pm 6, \ldots, \pm n= \pm 2m$ on 
the normal to the orbit.  
The refined stratification 
 $ \{ \tilde{S}_{\ub } : 
\ub  \in \tilde{{\cal B}}  \}$
 of $X$ this time has
as its strata the set $\tilde{S}_{(0)} = X^s$ of stable points, the set 
$\tilde{S}_{(T)}$ represented by
polynomials with exactly two distinct roots each with multiplicity $m=n/2$, 
for $2 \leq k \leq m$ the set $\tilde{S}_{(T,2k)}$ represented by polynomials
in the orbit of one of the form 
$$a_{m} t^{m} + a_{m+k}t^{m+k}  +  
a_{m+k+1}t^{m+k+1} + \cdots a_{2m} t^{2m}   $$
where $a_{m}$ and $a_{m+k}$
are nonzero, and finally the strata $S_{2j-n}$
(for $n/2 < j \leq n$) of the Morse stratification. 
\end{example}

\begin{example} \label{example5.3}
For a  more complicated example consider the action of $G=SL(3;\CC)$
and its maximal compact subgroup $K=SU(3)$ on $X=(\PP_2)^n$. A sequence
$(x_1,\ldots,x_n)$ of points in $\PP_2$ is semistable if and only if there is
no projective line $L$ in $\PP_2$ such that
\begin{equation} \label{line}
\frac{|\{j:x_j \in L \}|}{2} > \frac{n}{3}
\end{equation}
and no point $p \in \PP_2$ such that
\begin{equation} \label{point}
{|\{j:x_j = p \}|} > \frac{n}{3},
\end{equation}
and is stable if we can replace $>$ with $\geq$ in (\ref{line}) and (\ref{point}) (see
for example \cite{K2} (16.5)). The stratification $\{ S_\beta: \beta \in \calb \}$ can
be described as follows (see \cite{K2} Proposition 16.9). Any $(x_1, \ldots, x_n) \in (\PP_2)^n$
which is not semistable determines a unique flag
$$0 = M_0 \subset M_1 \ldots \subset M_s = \CC^3$$
in $\CC^3$ with $s=2$ or $3$, such that if $1 \leq i \leq s$ then
$$\frac{k_1}{m_1} > \cdots > \frac{k_s}{m_s}$$
where
$$k_i = |\{j: x_j \in M_i \backslash M_{i-1} \}| \mbox{ and } m_i = \dim(M_i/M_{i-1}),$$
and moreover if $m_i=2$ then those $x_j$ lying in $M_i \backslash M_{i-1}$ determine a semistable
element of $(\PP_1)^{k_i}$ after projection into $\PP(M_i/M_{i-1}) \cong \PP_1$. Then
$(x_1, \ldots, x_n)$ lies in the stratum labelled by the projection into $\mbox{Lie}(SU(3))$
of the vector
$$\beta = \left(\frac{k_1}{m_1}, \ldots , \frac{k_s}{m_s}\right) \in \mbox{Lie}(U(3))$$
in which each $k_i/m_i$ appears $m_i$ consecutive times. Thus $\calb$ is the projection
into the Lie algebra $\mbox{Lie}(SU(3))$ of
$$\{\left(\frac{n}{3},\frac{n}{3},\frac{n}{3}\right)\} \cup \{\left(\frac{k}{2},\frac{k}{2},n-k\right):\frac{2n}{3} < k \leq n\}
\cup \{ \left(k,\frac{n-k}{2}, \frac{n-k}{2}\right): \frac{n}{3} < k \leq n\}$$
$$ \cup \{(k_1,k_2,n-k_1 - k_2): n-k_1 - k_2 < k_2 < k_1 \leq n \}.$$
With a slight abuse of notation (identifying $\beta\in \mbox{Lie}(U(3))$ with its projection into  $\mbox{Lie}(SU(3))$)
we have
$$S_{(\frac{n}{3},\frac{n}{3},\frac{n}{3})} = X^{ss},$$
while if $2n/3 < k \leq n$ then $(x_1,\ldots,x_n)$ belongs to $S_{(k/2,k/2,n-k)}$
if and only if there is a line $L$ in $\PP_2$ containing exactly $k$ of the points $x_j$ and
at most $k/2$ of these points coincide anywhere on $L$. If $n/3 < k \leq n$ then
$(x_1,\ldots,x_n)$ belongs to $S_{(k,(n-k)/2,(n-k)/2)}$
if and only if exactly $k$ of the points $x_j$ coincide at some $p\in \PP_2$, and any line
through $L$ contains at most $(n-k)/2$ of the remaining points. Finally if
$n-k_1 - k_2 < k_2 < k_1 \leq n$ then $(x_1,\ldots,x_n)$ belongs to $S_{(k_1,k_2,n-k_1 - k_2)}$
if and only if there is a line $L$ in $\PP_2$ and a point $p \in L$ such that exactly $k_1$
of the points $x_j$ coincide at $p$ and exactly $k_2$ additional points lie on $L$.

The strata $S_\beta$ which need refining to obtain the stratification $\{ \tilde{S}_{\ub}:\ub \in
\tilde{\calb}\}$ are $S_{(k/2,k/2,n-k)}$ when $k$
is even and $S_{(k,(n-k)/2, (n-k)/2)}$ when $n-k$ is even, together with
$X^{ss}$ if $n$ is a multiple of $3$.  

If $k$ is even we have
$$S_{(\frac{k}{2},\frac{k}{2},n-k)} = \tilde{S}_{(\frac{k}{2},\frac{k}{2},n-k)} \sqcup
\tilde{S}_{(\frac{k}{2},\frac{k}{2},n-k,T_1)} \sqcup \tilde{S}_{(\frac{k}{2},\frac{k}{2},n-k,T_1,3)}$$
where $T_1 = \{ (t,t,t^{-2}): t \in \CC^*\}$. Here $(x_1,\ldots,x_n)$ belongs to the stratum $\tilde{S}_{(k/2,k/2,n-k)}$
if and only if there is a line $L$ in $\PP_2$ containing exactly $k$ of the points $x_j$ and
at most $k/2-1$ of these points coincide anywhere on $L$, while $(x_1,\ldots,x_n)$ belongs to 
$\tilde{S}_{(k/2,k/2,n-k,T_1)}$
if exactly $k/2$ of these points coincide somewhere on $L$ and the remaining $k/2$ coincide
somewhere else, and $(x_1,\ldots,x_n)$ belongs to $\tilde{S}_{(k/2,k/2,n-k,T_1,3)}$
if exactly $k/2$ of these points coincide somewhere on $L$ and the remaining $k/2$ on $L$
do not all coincide.

If $n-k$ is even we have 
$$S_{(k,\frac{n-k}{2},\frac{n-k}{2})} = \tilde{S}_{(k,\frac{n-k}{2},\frac{n-k}{2})} \sqcup
\tilde{S}_{(k,\frac{n-k}{2},\frac{n-k}{2},T_2)} \sqcup \tilde{S}_{(k,\frac{n-k}{2},\frac{n-k}{2},T_2,3)}$$
where $T_2 = \{ (t^{-2},t,t): t \in \CC^*\}$. Here $(x_1,\ldots,x_n)$ belongs to the stratum $\tilde{S}_{(k,(n-k)/2,(n-k)/2)}$
if and only if there is some point $p$ in $\PP_2$ where exactly $k$ of the points $x_j$ coincide, and
no line through $p$ contains at least $(n-k)/2$ of the remaining points $x_j$, 
while $(x_1,\ldots,x_n)$ belongs to 
$\tilde{S}_{(k,(n-k)/2,(n-k)/2,T_2)}$ if there are lines $L_1$ and $L_2$ meeting at $p$ and each containing
 exactly $(n-k)/2$ of the remaining points $x_j$, and $(x_1,\ldots,x_n)$ belongs to $\tilde{S}_{(k,(n-k)/2,(n-k)/2,T_2,3)}$
if there is a line $L$ through $p$ containing exactly $(n-k)/2$ of the remaining points $x_j$ but no other line
through $p$ contains at least $(n-k)/2$ of these remaining points.

Finally if $n$ is divisible by 3 then we have
$$X^{ss} = X^s \sqcup \tilde{S}_{(T)} \sqcup \{ \tilde{S}_{(T, \beta)}: \beta \in \calb_T \setminus \{ 0 \} \}
\sqcup \tilde{S}_{(T_1)} \sqcup \tilde{S}_{(T_1,3)} \sqcup \tilde{S}_{(T_1,-3)}$$
where $T$ is the standard maximal torus of $SU(3)$ and $\calb_T \setminus \{ 0 \}$ can be identified with
the set of four vectors
$$\{ \left(1,0,-1\right), \left(\frac{1}{2}, 0 , - \frac{1}{2}\right), \left(1, -\frac{1}{2}, -\frac{1}{2}
\right), \left(\frac{1}{2}, \frac{1}{2},-1\right) \}.$$
Here $(x_1,\ldots,x_n)$ belongs to $\tilde{S}_{(T)}$ if there are three points $p_1, p_2, p_3 \in \PP_2$
such that exactly $n/3$ of the points $x_j$ coincide at $p_i$ for $i=1,2,3$, while 
$(x_1,\ldots,x_n)$ belongs to $\tilde{S}_{(T_1)}$ if there is a point $p \in \PP_2$ and a line $L$ not
containing $p$ such that exactly $n/3$ of the points $x_j$ coincide at $p$ and the remaining $2n/3$
lie on $L$ with strictly fewer than $n/3$ coinciding anywhere on $L$. Also
$(x_1,\ldots,x_n)$ belongs to $\tilde{S}_{(T_1,3)}$ if there is a line $L$ in $\PP_2$ containing exactly
$2n/3$ of the points $x_j$ and strictly fewer than $n/3$ points coincide anywhere on $\PP_2$,
while $(x_1,\ldots,x_n)$ belongs to $\tilde{S}_{(T_1,-3)}$ if there is a point $p \in \PP_2$ where exactly
$n/3$ of the points $x_j$ coincide, and strictly fewer than $2n/3$ of the points $x_j$ lie on any line in $\PP_2$.

Finally, if $\beta \in \calb_T$ and  $(x_1,\ldots,x_n) \in \tilde{S}_{( T,\beta)}$ then there is some 
$p \in \PP_2$ and a line $L$ through $p$ such that exactly $n/3$ of the points $x_j$ coincide at
$p$ and exactly $2n/3$ lie on $L$, and
\begin{itemize}
\item[(a)] when $\beta = (1/2,0,-1/2)$ then $p$ and $L$ are unique;
\item[(b)] when $\beta = (1/2,1/2,-1)$ then $L$ is unique but there is another point $p'$ on $L$ where 
$n/3$ of the points $x_j$ coincide;
\item[(c)] when $\beta = (1,-1/2,-1/2)$ then $p$ is unique but there is another line $L'$ through $p$
containing exactly $2n/3$ of the points $x_j$;
\item[(d)] when $\beta = (1,0,-1)$ then there is another point $p' \in \PP_2 \setminus L$ where exactly $n/3$
of the points $x_j$ coincide, but the $n/3$ points $x_j$ lying on $L \setminus \{p\}$ do not all coincide.
\end{itemize}

\end{example}

\begin{example} \label{example5.4}
We can generalise Examples \ref{previo} and \ref{example5.3} by considering the action of 
$SL(m;\CC)$ on a product of Grassmannians
$$X = \prod_{j=1}^r \mbox{{\rm Grass}}(l_j, \CC^m)$$
where $\mbox{{\rm Grass}}(l,\CC^m)$ denotes the Grassmannian of $l$-dimensional linear
subspaces of $\CC^m$ and we linearise the action by using the Pl\"{u}cker embedding.
The Morse stratification $\{ S_\beta : \beta \in \calb \}$ is described in \cite{K2} $\S$16.3,
and the refinement $\tilde{S}_{\ub}: \ub \in \tilde{\calb} \}$ can be calculated by adapting the
methods used in \cite{INI} (especially \cite{INI} $\S$5) to refine the Yang--Mills stratification.

\end{example}

\section{Refinements when $G$ is abelian}
\renorm

The refinements we have considered so far of the Morse stratification $\{S_\beta: \beta \in \calb \}$ 
for $|\!|\mu |\!|^2$ are
unfortunately unlikely to be equivariantly perfect. When $G=T^c$ is abelian there is another way to 
refine this stratification which does lead to equivariantly perfect stratifications.

If $G$ is abelian then there is a partial desingularisation $X/\!/_\epsilon G$ of $X/\!/G$
obtained by perturbing the linearisation of the action of $G$ on $X$, or equivalently by replacing
the moment map $\mu :X \to \lieks$ by $\mu - \epsilon$ where $\epsilon \in \lieks$ is a generic
constant close to 0. Since $G$ is abelian the coadjoint action of $K=T$ on $\lieks = \liets$
is trivial and $\mu - \epsilon$ is an equivariant moment map for the action of $K$ on $X$.

Recall from $\S$1 that the Morse stratification $\{S_\beta: \beta \in \calb \}$ for $|\!|\mu |\!|^2$ 
is indexed by the closest points to 0 in the convex hulls of nonempty subsets of the set of weights
$\{\alpha_0,\ldots,\alpha_n\}$, and $x = [x_0:\ldots:x_n] \in X \subseteq \PP_n$ lies in the
stratum $S_\beta$ indexed by the closest point $\beta$ to 0 in the convex hull of
$$\{ \alpha_j: x_j \neq 0\}.$$
Similarly the Morse strata $\{S_\beta^\epsilon: \beta \in \calb^\epsilon \}$ for $|\!|\mu 
- \epsilon|\!|^2$ 
correspond to the closest points to $\epsilon$ in such convex hulls. More 
precisely $x = [x_0:\ldots:x_n] \in X $ lies in the
stratum $S_\beta^\epsilon$ indexed by the closest point $\beta$ to $0$ in the convex hull of
$\{ \alpha_j - \epsilon: x_j \neq 0\}$; then $\beta + \epsilon$ is the closest point to $\epsilon$
in the convex hull of 
$$\{ \alpha_j: x_j \neq 0\}.$$

If $\beta \in \calb$ does not lie in the convex hull of a proper subset of 
$\{ \alpha_j: \alpha_j.\beta = |\!|\beta |\!|^2 \}$, then for sufficiently
small $\epsilon$ we have
$$S_\beta = S_{\beta(\epsilon)}^\epsilon$$
for some small perturbation $\beta(\epsilon)$ of $\beta$, where $\beta(\epsilon) + \epsilon$ is the closest
point to $\epsilon$ of the convex hull of $\{ \alpha_j: \alpha_j.\beta = |\!|\beta |\!|^2 \}$. 
Otherwise for sufficiently small $\epsilon$ the stratum $S_\beta$ in the 
Morse stratification for $|\!|\mu |\!|^2$ is the union of $S_{\beta(\epsilon)}^\epsilon$
and other strata in the Morse stratification for $|\!|\mu 
- \epsilon|\!|^2$, which correspond to the closest points to $\epsilon$ of those
subsets of  $\{ \alpha_j: \alpha_j.\beta = |\!|\beta |\!|^2 \}$ whose closest point to 0
is $\beta$. Moreover if $\beta \in \calb^\epsilon$ and 
$\epsilon $ is chosen generically we can assume that $\beta$ does not lie in the
convex hull of any proper subset of $\{ \alpha_j: (\alpha_j- \epsilon).\beta = |\!|\beta |\!|^2 \}$.
This means that every point of
$$Z_\beta^\epsilon = \{ x \in X: x_j = 0 \mbox{ unless }(\alpha_j - \epsilon).\beta = |\!|\beta|\!|^2 \}$$
which is semistable for the action of 
the subgroup of $G$ whose Lie algebra is spanned by
$$\{\alpha_j - \epsilon - \beta: (\alpha_j - \epsilon) .\beta = |\!|\beta |\!|^2 \}$$
is also stable for the action of this subgroup. Moreover the set of semistable
points for the action of this subgroup is the same as the set 
$Z_\beta^{\epsilon,ss}$ of semistable points  of $Z^\epsilon_\beta$
under the action of the subgroup of $G$ whose Lie algebra is the
complexification of
the orthogonal complement to $\beta$ in $\liet$ (cf. Remark \ref{zbss}).

Thus we have a refinement of the Morse stratification for $|\!|\mu |\!|^2$ which is both
equivariantly perfect and has the property that all strata can be described in terms of the
sets of semistable points of nonsingular closed subvarieties $Z_\beta^\epsilon$ of $X$
under suitable linear actions of reductive subgroups of $G$ for which semistability coincides
with stability.

\begin{prop}
If $G=T^c$ is abelian and $\epsilon \in \lieks= \liets$ is chosen generically and sufficiently
close to 0 then the Morse stratification $\{S_\beta^\epsilon: \beta \in \calb^\epsilon \}$ for $|\!|\mu 
- \epsilon|\!|^2$ 
is an equivariantly perfect refinement of the Morse stratification $\{S_\beta: \beta \in \calb \}$ for $|\!|\mu 
|\!|^2$. In addition if $\beta \in \calb^\epsilon$ then we have
$$S_\beta^\epsilon = Y_\beta^{\epsilon,ss} = p_\beta^{-1}(Z_\beta^{\epsilon,ss})$$
and
$$H^*_G(S_\beta^\epsilon) = H^*_G(Z_\beta^{\epsilon,ss})$$
where
$$Z_\beta^{\epsilon} = \{ x \in X: x_j = 0 \mbox{ unless }(\alpha_j- \epsilon).\beta = |\!|\beta|\!|^2 \},$$
$$Y_\beta^\epsilon = \{ x \in X: x_j = 0 \mbox{ if }(\alpha_j - \epsilon).\beta < |\!|\beta|\!|^2 \mbox{ and 
$x_j \neq 0$ for some $j$ }$$
$$\mbox{ such that  }(\alpha_j - \epsilon).\beta = |\!|\beta|\!|^2\}$$
and
$$p_{\beta}(x) = \lim_{t \to \infty} \exp (-it\beta) x.$$
Furthermore $Z_\beta^{\epsilon,ss}$ is the set of semistable points for 
the action on $Z_\beta^\epsilon$ of the subgroup of $G$  whose Lie algebra is spanned by
$$\{\alpha_j - \epsilon - \beta: (\alpha_j - \epsilon) .\beta = |\!|\beta |\!|^2 \},$$
and every semistable point for this action is stable.
\end{prop} 

\begin{rem} When $G$ is abelian and $X^{ss}=X^s$ we can obtain the same result by perturbing the inner 
product on $\liek=\liet$ instead of perturbing the moment map from $\mu$ to $\mu - \epsilon$.
\end{rem}

\begin{rem}
Yet another possible procedure when $G = (\CC^*)^r$ is abelian is to use induction on $r$ to reduce to
the simple case when $G=\CC^*$.
\end{rem}

\begin{example} \label{example6.1}
Consider the action of the maximal torus
$$T^c = \{(t,t^{-1}):t \in \CC^* \} \cong \CC^*$$
of $G=SL(2;\CC)$ acting on $X=\PP_n$ or $X=(\PP_1)^n$ as in Examples \ref{previo} and \ref{example5.3}.
An element of $X$ represented by a sequence $(x_1,\ldots,x_n)$ of points of $\PP_1$ is semistable
(respectively stable) for $T^c$ if and only if at most $n/2$ (respectively strictly fewer than $n/2$) of the
points $x_j$ coincide at $0$ or $\infty$, and the Morse stratification $\{ S^T_\beta: \beta \in \calb^T \}$ for the
normsquare of the moment map $\mu_T: X \to \liet$ is indexed by
$$\{ 0 \} \cup \{ 2j-n: 0 \leq j \leq n \}$$
with $(x_1,\ldots,x_n)$ representing a point of $S^T_{2j-n}$ when exactly $j$ of the points $x_i$ coincide at
0 if $j>n/2$, and exactly $n-j$ of the points $x_i$ coincide at $\infty$ if $j<n/2$. When $n$ is odd this
stratification is unchanged if we perturb the linearisation slightly. When $n$ is even the only change is that
$S_0^T$ becomes two strata: the subset represented by $(x_1,\ldots,x_n)$ with exactly $n/2$ of the points
coinciding at 0 (or at $\infty$, depending on the sign of the perturbation), and its complement in $S_0^T$.
This contrasts with the refinement $\{ \tilde{S}^T_{\ub}: \ub \in \tilde{\calb}^T \}$ in which $S_0^T$ is subdivided into
three strata when $X=(\PP_1)^n$ and more when $X=\PP_n$ (cf. Examples \ref{previo} and \ref{example5.3}.
\end{example}

\begin{example}
When the maximal torus $T^c \cong (\CC^*)^2$ of $SL(3;\CC)$ acts on $X=(\PP_2)^n$ then the Morse stratification
$\{S^T_\beta: \beta \in \calb^T \}$ for $|\!| \mu_T |\!|^2$ and its refinement $\{\tilde{S}^T_{\ub}: \ub \in \tilde{\calb}^T \}$
defined as in $\S$5 can be described in a way closely analogous to the stratifications
$\{S_\beta: \beta \in \calb \}$ and  $\{\tilde{S}_{\ub}: {\ub} \in \tilde{\calb} \}$ associated to the action of $G=SL(3;\CC)$
which were described in Example \ref{example5.4}; the difference is that the points $p$ and lines $L$ which
appear in the description are $[1:0:0]$, $[0:1:0]$ or $[0:0:1]$ or one of the lines joining these points. A generic
perturbation of the linearisation provides a coarser refinement of 
$\{S^T_\beta: \beta \in \calb^T \}$ than the stratification $\{\tilde{S}^T_{\ub}: \ub \in \tilde{\calb}^T \}$.
If $k$ or $n-k$ is even then the strata $S^T_\beta$ with indices $\beta$ in the Weyl group orbit of
$(k/2,k/2,n-k)$ or $(k,(n-k)/2, (n-k)/2)$ decompose as the disjoint union of two refined strata instead of the three
in the refinement  $\{\tilde{S}^T_{\ub}: \ub \in \tilde{\calb}^T \}$, in a way which is directly analogous to the 
decomposition of $S_0^T$ in Example \ref{example6.1}. Similarly if $n$ is divisible by 3 then $S^T_0$
decomposes into a smaller number of strata than in the stratification
 $\{\tilde{S}^T_{\ub}: \ub \in \tilde{\calb}^T \}$; in particular the analogues of the strata $\tilde{S}_{(T)}$ and
$\tilde{S}_{(T_1)}$ described in Example \ref{example5.4} are amalgamated with higher dimensional strata.

\end{example}

\section{Reduction to a maximal torus}

\renorm

The Morse stratification $\{ S_\beta: \beta \in \calb \}$ of the normsquare of the moment map
is useful for studying the cohomology of the GIT quotient $X/\! /G$ because it is equivariantly 
perfect over the rationals; that is, if $\lambda(\beta)$ is the real codimension of $S_\beta$ in $X$ and if
$U_\beta$ is the open subset
$$U_\beta = S_\beta \cup \bigcup_{|\!| \beta' |\!| < |\!| \beta |\!| } S_{\beta'}$$
which contains $S_\beta$ as a closed subset, then the Thom--Gysin long exact sequences
$$\ldots \to H_G^{*-\lambda(\beta)}(S_\beta) \stackrel{TG_\beta}{\to} H_G^*(U_\beta) \to H^*_G(U_\beta \setminus S_\beta)
\to \ldots$$
break up into short exact sequences
$$0 \to H_G^{*-\lambda(\beta)}(S_\beta) \stackrel{TG_\beta}{\to} H_G^*(U_\beta) \to H^*_G(U_\beta \setminus S_\beta)
\to 0$$
of equivariant cohomology with rational coefficients.

\begin{rem} \label{restricb}
This happens because the composition of the Thom--Gysin map $TG_\beta$ with the restriction map
$H^*_G(U_\beta) \to H^*_G(S_\beta)$ is given by multiplication by the equivariant Euler class $e_\beta$
of the normal bundle to $S_\beta$ in $U_\beta$, and it follows from a criterion of Atiyah and Bott
\cite{AB} $\S$13 that $e_\beta$ is not a zero-divisor in $H^*_G(S_\beta)$. Note that the restriction
maps from $H^*_G(X)$ to $H^*_G(X^{ss})$ and also to $H^*_G(U_\beta)$ for $\beta \in \calb \setminus \{ 0 \}$
are compositions of restriction maps from $H^*_G(U_{\beta'})$ to $H^*_G(U_{\beta'}\setminus S_{\beta'})$ 
which are all surjective. In particular the cohomology ring $H^*_G(X^{ss}$ (which in the good
case when $X^{ss} = X^s$ is isomorphic to the rational cohomology ring of the 
geometric invariant theoretic quotient $X/\!/G$ or equivalently the Marsden-Weinstein reduction $\mu^{-1}(0)/K$)
is isomorphic to the quotient of $H^*_G(X)$ by the kernel of the restriction map $\rho:H^*_G(X) \to H^*_G(X^{ss})$.
\end{rem}

Unfortunately the refinement $\{ \tilde{S}_{\ub}: \ub \in \tilde{\calb} \}$ described in $\S$5 
of $\{ S_\beta: \beta \in \calb \}$ is not in
general equivariantly perfect. It is still the case that we have Thom--Gysin long exact
sequences for equivariant cohomology, and the kernels of the restriction maps
from $H^*_G(X)$ to $H_G^*(X^{ss})$ and to $H^*_G(X^s)$ can be described in terms of the
images of the associated Thom--Gysin maps, but the description is not as clean as in
the equivariantly perfect case, and the restriction map to $H^*_G(X^s)$ is not necessarily
surjective (see Example \ref{lastex} below). 

However when $G$ is abelian the alternative refinement
$\{ S_\beta^\epsilon: \beta \in \calb^\epsilon \}$ described in $\S$6 is equivariantly perfect. We
can exploit this fact even when $G$ is not abelian by making use of the close relationship
between the restriction maps
$$\rho:H^*_G(X) \to H^*_G(X^{ss})$$
and
$$\rho_T: H^*_{T}(X) \to H^*_{T}(X^{ss,T})$$
where $T^c$ is a maximal torus of $G$ (we can assume that it is the complexification of
$T=T^c \cap K$ which is a maximal torus of $K$) and $X^{ss,T}$ is the set of $T^c$-semistable
points of $X$. Therefore in this final section we shall describe this relationship.

The kernel of $\rho$ can be described in several different ways. The first follows 
immediately from the fact that the stratification $\{ S_\beta: \beta \in \calb \}$ is equivariantly perfect. 

\begin{lem} \label{AA}
For $\beta \in \calb \setminus \{ 0 \}$ let $\widetilde{TG}_\beta:H^{*-\lambda(\beta)}_G(S_\beta) \to
H^*_G(X)$ be any lift to $X$ of the Thom--Gysin map
${TG}_\beta:H^{*-\lambda(\beta)}_G(S_\beta) \to
H^*_G(U_\beta)$. Then
$$\ker \rho = \bigoplus_{\beta \in \calb \setminus \{ 0 \}} {\rm im} \widetilde{TG}_\beta.$$
\end{lem}

\begin{rem}
Note that such lifts always exist by Remark \ref{restricb}.
\end{rem}

Another closely related description is given in \cite{E,EK,K}.

\begin{lem} \label{BB}
Let $\calr$ be a subset of $\ker \rho$ such that for every $\beta \in \calb \setminus \{ 0 \}$ and
every $\eta \in H^*_G(S_\beta)$ there exists $\zeta \in \calr$ such that
$$ \zeta|_{S_\gamma} = 0 \mbox{ unless } |\!|\gamma|\!| \geq |\!|\beta |\!|$$
and
$$\zeta|_{S_\beta} = \eta e_\beta$$
where $e_\beta$ is the equivariant Euler class of the normal bundle to $S_\beta$. Then $\calr$ spans
$\ker \rho$.
\end{lem}

In the good case when $X^{ss} = X^s$ a rather different description follows from the formulas for the
intersection pairings in $H^*(X/\!/G)$ given in \cite{JK2} (see also \cite{GK,M,M1,WI}).

\newcommand{\res}{\mbox{res}}

\begin{lem} \label{CC}
If $X^{ss} = X^s$ and $\eta \in H^j_G(X)$ then $\eta \in \ker \rho$ if and only if
$$ \res \Biggl ( 
\cald^2 
 \sum_{F \in \calf} e^{ \mu_F}\int_F  \frac{( \eta \zeta)|_F  
 }{e_F  }  \Biggr ) = 0$$
for all $\zeta \in H^k_G(X)$ with $j+k = \dim_\RR (X/\!/G)$. 
Here $\calf$ is the set of connected components $F$ of the fixed point set of 
$T$ in $X$ 
 and $e_F $ 
$\in H^*_T(X)$ is the $T$-equivariant Euler class of the normal 
bundle to $F$ in $X$. Also $\mu_F:\liet \to \RR$ denotes the
linear function on $\liet$ given by the constant value in $\liets$
taken by the moment map $\mu$  on  
$F \in \calf$, and the polynomial
$\cald: \liet \to \RR$ is defined by 
$$\cald = \prod_{\gamma > 0 } \gamma,$$
where $\gamma$ runs 
over the positive roots of $K$. 
\end{lem}

\noindent{{\bf Proof}:} This follows immediately from Poincar\'{e} duality for  the orbifold 
$X/\!/G$ and the surjectivity of 
$$\rho:H^*_G(X) \to H^*_G(X^{ss}) \cong H^*(X/\!/G)$$
together with \cite{JK2} Theorem
8.1, which tells us that, up to multiplication by a nonzero constant,
the formula in Lemma \ref{CC} is the intersection pairing 
$$ \int_{X/\!/G} \rho(\eta) \rho(\zeta)$$
of $\rho(\eta)$ and $\rho(\zeta)$ in $X/\!/G$.

\begin{rem} The multivariable residue $\res$ which appears in Lemma \ref{CC} above 
is a linear map defined on a class of  meromorphic differential forms on $\liet\otimes\CC$.
In order to apply it to the individual terms in the residue
formula it is necessary to make some choices which do not affect the residue
of the whole sum. Once the choices have been made, many
of the terms in the sum contribute zero and the formula
 can be rewritten as
a sum over a subset $\calf_+$ of the set $\calf $ of components of the
fixed point set $X^T$, consisting of those $F \in \calf$ on which the constant
value taken by $\mu_T$ lies in a certain cone with its vertex at 0. 
When $\dim T =1$ and $\liet$ is identified
with $\RR$, we can take
$$\calf_+ = \{ F \in \calf: \mu_T(F) > 0 \}$$
and if we replace $\calf$ with $\calf_+$ in the formula in Lemma \ref{CC} 
then we can interpret $\res$ as the usual residue at 0 of a rational function
in one variable.
\end{rem}

\begin{rem} \label{DD}
When $G$ is abelian we have analogues of Lemmas \ref{AA} and \ref{BB} for the 
refinement 
$\{ S_\beta^\epsilon: \beta \in \calb^\epsilon \}$ of the stratification
$\{S_\beta: \beta \in \calb \}$ obtained by perturbing the
moment map $\mu$ to $\mu - \epsilon$ (see $\S$6). We also have the
following description of $\ker \rho$ due to Tolman and Weitsman \cite{TW}.
We can use these results even when $G$ is not abelian via Lemma \ref{FF}
below. 

\end{rem}

\begin{lem} \label{EE}
When $G=T^c$ is abelian and $X^{ss}=X^s$ then
$$\ker \rho = \sum_{\xi \in \liet} \{\eta \in H^*_T(X): \eta|_{F \cap X_\xi} = 0 \mbox{ for all } F \in \calf \}
$$
where
$$X_\xi = \{ x \in X: \mu(x)(\xi) \leq 0 \}$$
if $\xi \in \liet$.
\end{lem}

Recall that the Weyl group $W$ of $K$ acts faithfully on the Lie algebra $\liet$ of $T$
and is generated by reflections. For any $w \in W$ we denote by $(-1)^w$ the determinant
of $w$ regarded as an automorphism of $\liet$. If $W$ acts on a module $M$ then
$$M^W = \{m\in M: wm=m \mbox{ for all } w \in W \}$$
is the set of $W$-invariants in $M$ and
$$M^{antiW} = \{ m \in M: wm=(-1)^w m \mbox{ for all } w \in W \}$$
is the set of anti-invariants for the action of $W$ on $M$. There is a natural identification of
$H^*_G(X)$ with $[H^*_T(X)]^W$, which in the case when $X$ is a point becomes a natural
identification of $H^*(BG)$ with $[H^*(BT)]^W$ where $H^*(BT)$ is the $\QQ$-algebra
of polynomial functions on $\liet$. There is (up to sign) one fundamental anti-invariant
in $H^*(BT)$, which is the product $\cald$ of the positive roots of $K$, and we have
the following wellknown facts (see for example \cite{ES} Lemma 1.2).

\begin{lem} \label{wantiw} (i) $[H^*(BT)]^{antiW}$ is a free $[H^*(BT)]^W$-module
of rank one generated by $\cald$.

\noindent (ii) $[H^*(BT)]^{antiW}$ is a direct summand of $H^*(BT)$ as an
 $[H^*(BT)]^W$-module, with splitting given by

$$\eta \mapsto \frac{1}{|W|} \sum_{w \in W} (-1)^w w\eta.$$
\end{lem}

Note that $H^*_T(X)$ has analogous properties, since it is isomorphic to
$H^*(BT) \otimes H^*(X)$ as an $H^*(BT)$-module (see \cite{K2} Proposition
5.8), though not as a ring.

\begin{lem} \label{FF}
Suppose that $X^{ss} = X^s$. If $\eta \in H^*_G(X) \cong [H^*_T(X)]^W$ then
the following are equivalent:

(i) $\eta \in \ker \rho;$

(ii) $\eta \cald \in \ker \rho_T;$

(iii) $\eta \cald^2 \in \ker \rho_T.$

\noindent Moreover multiplication by $\cald$ induces a bijection
\begin{equation} \label{bijrrt} \ker \rho \to \ker \rho_T \cap [H^*_T(X)]^{antiW} \end{equation}
with inverse
$$ \eta \mapsto \frac{1}{\cald |W|} \sum_{w \in W} (-1)^w w\eta.$$
\end{lem}

\begin{rem}
The corresponding results are true when $X$ is a compact symplectic
manifold with a Hamiltonian action of a compact group $K$, provided that
0 is a regular value of the moment map.
\end{rem}

\begin{rem} \label{rem711}
Martin \cite{M,M1} gave a direct proof of the equivalence $(i) \Leftrightarrow (iii)$ when
$X^{ss} = X^s$, while Guillemin and Kalkman \cite{GK} observed that it follows immediately
from Lemma \ref{CC}. The equivalence $(i)\Leftrightarrow (ii)$ and the bijection (\ref{bijrrt})
goes back in essence at least to Ellingsrud and Str\o mme \cite{ES} $\S$4. Direct proofs of the
equivalences $(i) \Leftrightarrow (ii)$ and $(ii) \Leftrightarrow (iii)$ and the bijection (\ref{bijrrt})
are given below;
the assumption that $X^{ss} = X^s$ is not needed for the equivalence $(i) \Leftrightarrow (ii)$
or for the existence of the bijection (\ref{bijrrt}), and in fact the most convenient assumption for
the proof below of the equivalence $(ii) \Leftrightarrow (iii)$ is not that $X^{ss} = X^s$ but rather
the corresponding assumption $X^{ss,T}=X^{s,T}$ for the action of the maximal torus. Thus
Lemma \ref{FF} is true when either $X^{ss} = X^s$ or $X^{ss,T} = X^{s,T}$, and if (iii) is omitted
then it is true without either of these hypotheses.

Note also that the direct proofs of the equivalences 
$(i) \Leftrightarrow (iii)$ and $(ii) \Leftrightarrow (iii)$
do not require $X$ to be compact; it suffices that $\mu^{-1}(0)$
or equivalently $X/\!/G$ should be compact, and for example
they can be applied to the study of moduli spaces of
bundles over a compact Riemann surface (cf. \cite{JK2}).
\end{rem}

Before giving a proof of Lemma \ref{FF} we shall relate the stratifications $\{S_\beta:\beta \in \calb \}$
and $\{ S^T_\beta : \beta \in \calb_T \}$ and the associated Thom--Gysin maps $TG_\beta$ and
$TG^T_\beta$ for the actions of $G$ and its maximal torus $T^c$ on $X$. Recall from $\S$1 that the
indexing set $\calb_T$ consists of the closest points to 0 of the convex hulls in $\liets \cong \liet$
of the weights $\alpha_0, \ldots, \alpha_n$ for the linear action of $T^c$ on $X \subseteq \PP_n$, and
that
$$\calb = \calb_T \cap \liet_+ \mbox{ and } \calb_T = \{ w\beta: \beta \in \calb, w \in W \}$$
where $\liet_+$ is a positive Weyl chamber in $\liet$. If $\beta \in \calb$ then in the notation of $\S$1
we have
$$S_\beta = G Y_\beta^{ss} \mbox{ and } S_\beta^T = Y_\beta^{ss,T}$$
where $Y_\beta^{ss,T} = p_\beta^{-1}(Z_\beta^{ss,T})$ and $Z_\beta^{ss,T}$ is the set of semistable
points for an appropriate linearisation of the action of $T^c$ on $Z_\beta$. There are induced
isomorphisms
$$H^*_G(S_\beta) \cong H^*_{\stab(\beta)} (Z_\beta^{ss}) \mbox{ and } H^*_T(S_\beta^T) \cong H^*_T(Z_\beta^{ss,T})$$
with surjective restriction maps
$$ H^*_{\stab(\beta)}(Z_\beta) \cong [H_T^*(Z_\beta)]^{W_\beta} \to H^*_{\stab(\beta)}(Z_\beta^{ss})$$
and
$$H^*_T(Z_\beta) \to H^*_T(Z_\beta^{ss,T})$$
where $W_\beta$ is the Weyl group of $\stab(\beta)$. Notice that the boundary
$\bar{S}^T_\beta \setminus S^T_\beta$ of $S^T_\beta$ is contained in 
$$\bigcup_{\beta' \in \calb_T: |\!| \beta' |\!| > |\!| \beta |\!|} S^T_{\beta'} \subseteq 
\bigcup_{\beta' \in \calb: |\!| \beta' |\!| > |\!| \beta |\!|} S_{\beta'}.$$
Let $\cald_\beta \in H^*(BT)$ be the product of the positive roots of $\stab(\beta)$, and if $\beta \in \calb$
let $\widetilde{TG}^T_\beta$ be any lift to $H^*_T(X)$ of the Thom--Gysin map
${TG}^T_\beta$ associated to the inclusion of $S^T_\beta$ in 
$$U^T_\beta = S^T_\beta \cup \bigcup_{\beta' \in \calb_T: |\!| \beta' |\!| < |\!| \beta |\!|} S^T_{\beta'}.$$
If $\eta \in H^*_T(Z_\beta)$ then $\eta$ restricts to an element of $H^*_T(Z_\beta^{ss,T}) \cong 
H^*_T(Y_\beta^{ss,T}) = H^*_T(S_\beta^T)$ and 
$$\sum_{w \in W} (-1)^w w\left( \cald_\beta \widetilde{TG}_\beta^T(\eta) \right)$$
is a $W$-anti-invariant element of $H^*_T(X)$. Thus we have a well defined element
$$\frac{1}{\cald} \sum_{w \in W} (-1)^w w\left( \cald_\beta \widetilde{TG}_\beta^T(\eta) \right)$$
of $[H^*_T(X)]^W \cong H^*_G(X)$ (cf. Lemma \ref{wantiw}), whose restriction to
$$U_\beta = S_\beta \cup \bigcup_{\beta' \in \calb: |\!| \beta' |\!| < |\!| \beta |\!|} S_{\beta'}$$
is independent of the choice of lift $\widetilde{TG}^T_\beta$ of the Thom--Gysin map
${TG}^T_\beta$ and thus can be expressed as
$$\frac{1}{\cald} \sum_{w \in W} (-1)^w w\left( \cald_\beta {TG}_\beta^T(\eta) \right)
 = \frac{1}{\cald} \sum_{w \in W} (-1)^w \cald_{w\beta} {TG}_{w\beta}^T(w\eta).$$

\begin{lem} \label{GG}
If $\beta \in \calb$ and $\eta \in H^*_{\stab(\beta)}(Z_\beta) \cong [H^*_T(Z_\beta)]^{W_\beta}$
then
$$TG_\beta (\eta) = \frac{1}{|W|\cald} \sum_{w \in W} (-1)^w w\left( \cald_\beta {TG}_\beta^T(\eta) \right)
 = \frac{1}{|W|\cald} \sum_{w \in W} (-1)^w \cald_{w\beta} {TG}_{w\beta}^T(w\eta).$$
\end{lem}

\noindent{\bf Proof}: Note that $\eta \in H^*_{\stab(\beta)}(Z_\beta)$ represents an element of
$H^*_{\stab(\beta)}(Z_\beta^{ss}) \cong H^*_G(S_\beta)$
 by restriction from $Z_\beta$ to
$Z_\beta^{ss}$. Since the equivariant Euler class $e_\beta$ of the normal bundle $\caln_\beta$
to $S_\beta$ is not a zero-divisor in $H^*_G(S_\beta)$, it suffices to show that 
$$
\frac{1}{|W|\cald} \sum_{w \in W} (-1)^w w\left( \cald_\beta {TG}_\beta^T(\eta) \right)
 = \frac{1}{|W|\cald} \sum_{w \in W} (-1)^w \cald_{w\beta} {TG}_{w\beta}^T(w\eta)$$
restricts to 0 on 
$$\bigcup_{\beta' \neq \beta, |\!| \beta' |\!| \leq |\!| \beta |\!|} S_{\beta'}$$
and restricts to $\eta e_\beta$ on $S_\beta$. 
Since 
$H^*_{\stab(\beta)}(Z_\beta^{ss}) \cong H^*_G(S_\beta)$
it suffices to check that the restriction to $Z_\beta^{ss}$ is $\eta e_\beta$.
But 
$$\overline{S^T_{w\beta}} = w(\overline{S^T_\beta}) \subseteq S_\beta \cup \bigcup_{|\!|\beta'|\!|
> |\!| \beta |\!|} S_{\beta'}$$
so $TG^T_{w\beta}(w \eta)$ restricts to 0 on
$$\bigcup_{\beta' \neq \beta, |\!| \beta' |\!| \leq |\!| \beta |\!|} S_{\beta'}$$
as required. Also the composition of $TG^T_\beta$ with restriction to $S_\beta^T$
is multiplication by the equivariant Euler class $e^T_\beta$ of the normal
bundle $\caln^T_\beta$ to $S^T_\beta$, and
$$Z_\beta^{ss} \subseteq Z_\beta^{ss,T} \subseteq S_\beta^T.$$
Hence the restriction of $\cald_\beta TG^T_\beta(\eta)$ to $Z_\beta^{ss}$ is
$\cald_\beta \eta e^T_\beta$. Now
$$S_\beta \cong G \times_{P_\beta} Y^{ss}_\beta \mbox{ and } S_\beta^T = Y_\beta^{ss,T}$$
where $Y_\beta^{ss}$ is an open subset of $Y_\beta^{ss,T}$, so their normal bundles
$\caln_\beta$ and $\caln^T_\beta$ are related by
$$\caln^T_\beta|_{Y_\beta^{ss}} \cong \frak{g}/{\frak{p}_\beta} \oplus \caln_\beta|_{Y_\beta^{ss}}$$
where $\frak{g}$ and $\frak{p}_\beta$ are the Lie algebras of $G$ and $P_\beta$. Therefore
on restriction to $Y_\beta^{ss}$ (or to $Z_\beta^{ss}$) we have
$$e_\beta^T = \frac{\cald e_\beta}{\cald_\beta} $$
so the restriction of $\cald_\beta TG^T_\beta(\eta)$ to $H^*_{\stab(\beta)}(Z_\beta^{ss})$ is
$$\cald_\beta \eta e_\beta^T = \cald \eta e_\beta.$$
Hence the restriction to $H^*_G(S_\beta)$ of
$$
\frac{1}{|W|\cald} \sum_{w \in W} (-1)^w w\left( \cald_\beta {TG}_\beta^T(\eta) \right)
 = \frac{1}{|W|\cald} \sum_{w \in W} (-1)^w \cald_{w\beta} {TG}_{w\beta}^T(w\eta)$$
is 
$$
\frac{1}{|W|\cald} \sum_{w \in W} (-1)^w w\left( \cald \eta e_\beta) \right) = \eta e_\beta$$
as required, since $e_\beta$ and $\eta \in H_G^*(S_\beta)$ are $W$-invariant.

\bigskip

\noindent{\bf Proof of Lemma \ref{FF}}: First recall from Lemma \ref{wantiw} that the map
$$p:[H^*_T(X)]^{antiW} \to [H^*_T(X)]^W \cong H^*_G(X)$$
defined by
$$p(\zeta) = \frac{1}{|W|\cald} \sum_{w \in W} (-1)^w w\zeta$$
is a bijection whose inverse is given by multiplication by $\cald$, since if $\zeta$ is
$W$-invariant then $p(\cald \zeta) = \zeta$ and if $\zeta$ is anti-invariant then
$\cald p(\zeta) = \zeta$ (cf. \cite{ES} (4.3)). Suppose that $\eta \in H^*_G(X)$. If
$\eta \in \ker \rho$, then by Lemma \ref{AA} we can write
$$\eta = \sum_{\beta \in \calb \setminus \{ 0 \}} \widetilde{TG}_\beta(\eta_\beta)$$
for some $\eta_\beta \in H^*_{\stab(\beta)}(Z_\beta)$ representing an element
of $H^*_G(S_\beta) \cong H^*_{\stab(\beta)}(Z_\beta^{ss})$, where 
$\widetilde{TG}_\beta:H^{*-\lambda(\beta)}_G(S_\beta) \to
H^*_G(X)$ is any lift to $X$ of the Thom--Gysin map
${TG}_\beta:H^{*-\lambda(\beta)}_G(S_\beta) \to
H^*_G(U_\beta)$. By Lemma \ref{GG} we can choose $\widetilde{TG}_\beta$ so that
$$\widetilde{TG}_\beta(\eta_\beta) = 
\frac{1}{|W|\cald} \sum_{w \in W} (-1)^w \cald_{w\beta} {TG}_{w\beta}^T(w\eta_\beta)$$
where $\widetilde{TG}^T_\beta$ is a lift of the Thom--Gysin map $TG^T_\beta$ and the
choice of these lifts $\widetilde{TG}^T_\beta$ respects the action of $W$. Then
$$\cald \eta = \sum_{\beta \in \calb \setminus \{ 0 \}} 
\frac{1}{|W|} \sum_{w \in W} (-1)^w \cald_{w\beta} {TG}_{w\beta}^T(w\eta_\beta)
\in \bigoplus_{\beta \in \calb_T \setminus \{ 0 \}} {\rm im} \widetilde{TG}^T_\beta = \ker \rho_T.$$
Conversely, suppose that $\cald \eta \in \ker \rho_T$. Then by Lemma \ref{AA}
\begin{equation} \label{up}
\cald \eta = \sum_{\beta \in \calb_T \setminus \{ 0 \}} \widetilde{TG}^T_\beta(\eta_\beta)
\end{equation}
for some $\eta_\beta \in H^*_T(Z_\beta)$ representing an element of 
$H^*_T(S^T_\beta) \cong H^*_T(Z_\beta^{ss,T})$, where 
$\widetilde{TG}^T_\beta$ is any lift to
$H^*_G(X)$ of the Thom--Gysin map
${TG}_\beta^T$. But $\eta$ is $W$-invariant, so $\cald \eta$ is anti-invariant, and so
$$\cald \eta = \frac{1}{|W|} \sum_{w_1 \in W} (-1)^{w_1} w_1(\cald \eta)
= \frac{1}{|W|} \sum_{\beta \in \calb_T \setminus \{ 0 \}} 
\sum{w_1 \in W} (-1)^{w_1} \widetilde{TG}^T_{w_1 \beta}(w_1 \eta_\beta)$$
if the lifts $\widetilde{TG}^T_\beta$ are chosen to respect the action of $W$.
We have
$$\calb_T = \{w\beta: w \in W, \beta \in \calb \}$$
and if $w \in W$ and $\beta \in \calb$ then $w\beta = \beta$ if and only if 
$w \in W_\beta$. Thus
$$\cald \eta 
= \frac{1}{|W|} \sum_{\beta \in \calb_T \setminus \{ 0 \}} \sum_{w_2 \in W} 
\sum{w_1 \in W} \frac{(-1)^{w_1}}{|W_\beta|}  \widetilde{TG}^T_{w_1 w_2 \beta}(w_1 \eta_{w_2\beta})$$
$$= \frac{1}{|W|} \sum_{\beta \in \calb_T \setminus \{ 0 \}} \sum_{w_2 \in W} 
\sum{w \in W} \frac{(-1)^w (-1)^{w_2}}{|W_\beta|}  \widetilde{TG}^T_{w \beta}(w w_2^{-1} \eta_{w_2\beta})$$
$$=  \frac{1}{|W|} \sum_{\beta \in \calb_T \setminus \{ 0 \}} \sum_{w \in W} (-1)^w w \left(
 \widetilde{TG}^T_{ \beta}\left( \sum{w_2 \in W} \frac{(-1)^{w_2}}{|W_\beta|} w_2^{-1} \eta_{w_2\beta}\right)\right).
$$
If $\tilde{w} \in W_\beta$ then
$$\tilde{w} \left( \sum_{w_2 \in W} \frac{(-1)^{w_2}}{|W_\beta|} w_2^{-1}\eta_{w_2 \beta} \right)
= \sum_{w_2 \in W} \frac{(-1)^{w_2}}{|W_\beta|} \tilde{w} w_2^{-1}\eta_{w_2 \beta}$$
$$ =  \sum_{w_3 \in W} \frac{(-1)^{w_3}(-1)^{\tilde{w}}}{|W_\beta|} w_3^{-1}\eta_{w_3\tilde{w} \beta}
=  (-1)^{\tilde{w}}  \sum_{w_3 \in W} \frac{(-1)^{w_3}}{|W_\beta|} w_3^{-1}\eta_{w_3 \beta}.$$
Since elements of $H^*_T(Z_\beta)$ which are anti-invariant under the action of $W_\beta$
are multiples of $\cald_\beta$, it follows that
$$\sum_{w_2 \in W} \frac{(-1)^{w_2}}{|W_\beta|} w_2^{-1}\eta_{w_2 \beta} = \cald_\beta \zeta_\beta$$
for some $\zeta_\beta \in [H^*_T(Z_\beta)]^{W_\beta} \cong H^*_{\stab(\beta)}(Z_\beta)$ and hence
$$\cald \eta =  \frac{1}{|W|} \sum_{\beta \in \calb_T \setminus \{ 0 \}} \sum_{w \in W} 
(-1)^w w \left( \cald_\beta \widetilde{TG}^T_\beta(\zeta_\beta)\right)$$
so that
$$\eta 
\in \bigoplus_{\beta \in \calb \setminus \{ 0 \}} {\rm im} \widetilde{TG}_\beta = \ker \rho$$
by Lemmas \ref{AA} and \ref{GG}. This proves the equivalence $(i) \Leftrightarrow (ii)$, and the
same argument shows that the bijection
$$H^*_G(X) \to [H^*_T(X)]^{antiW}$$
given by multiplication by $\cald$ restricts to a bijection from $\ker \rho$ to $\ker \rho_T \cap [H^*_T(X)]^{antiW}$.

The observation that the  equivalence $(i) \Leftrightarrow (iii)$ follows directly from Lemma \ref{CC} 
when $X^{ss} = X^s$ now completes the proof of Lemma \ref{FF}. However it is also easy to show
directly that if $X^{ss,T} = X^{s,T}$ then $(ii) \Leftrightarrow (iii)$ for any $\eta \in H^*_G(X) \cong [H^*_T(X)]^W$;
that is, $\cald \eta \in \ker \rho_T$ if and only if $\cald^2 \eta \in \ker \rho_T$. This follows from
Lemma \ref{CC} applied to the action of $T^c$, together with Poincar\'{e} duality on $X/\!/T^c$
and the surjectivity of $\rho_T$, as
$$\cald \eta \in \ker \rho_T \Leftrightarrow \int_{X/\!/T^c} \rho_T(\cald \eta) \rho_T(\zeta) = 0
\mbox{ for all }\zeta \in H^*_T(X).$$
Since $\cald \eta$ is anti-invariant, this holds for all $\zeta \in H^*_T(X)$ if and only if it holds
for all anti-invariant $\zeta$; that is, for all $\zeta$ of the form $\cald \xi$ where $\xi \in [H^*_T(X)]^W$.
Thus 
$$\cald \eta \in \ker \rho_T \Leftrightarrow 0 = \int_{X/\!/T^c} \rho_T(\cald \eta) \rho_T(\cald \xi) = 
 \int_{X/\!/T^c} \rho_T(\cald^2 \eta \xi) $$
for all invariant $\xi \in H^*_T(X)$, or equivalently (since $\cald^2 \eta$ is invariant) for all
$\xi \in H^*_T(X)$. But by Poincar\'{e} duality again we have 
$$\int_{X/\!/T^c} \rho_T(\cald^2 \eta \xi) =0$$
for all $\xi \in H^*_T(X)$ if and only if $\rho_T(\cald^2 \eta) = 0$, as required.

\begin{rem} 
These ideas are used in \cite{EK} to obtain a complete set of relations between the
standard generators of the moduli space $\mnd$ of stable holomorphic vector bundles
of rank $n$ and degree $d$ over a fixed compact Riemann surface of genus $g \geq 2$
when $n$ and $d$ are coprime. There the r\^{o}le of $X/\!/G$ is played by the moduli
space $\mnd$, and the r\^{o}le of $X/\!/T^c$ is played by the corresponding moduli
space of parabolic bundles where the parabolic structure is associated to a full flag.
The generic perturbation of the linearisation used in $\S$6 to obtain a refined stratification
is played by a generic perturbation of the parabolic weights.
\end{rem}

\begin{example}
Suppose that $X=\PP_n$ and that as usual the maximal torus $T^c$ of $G$ acts with
weights $\alpha_0, \ldots, \alpha_n$. In this case the closure $\bar{S}^T_\beta = \overline{Y_\beta}$
of any $T^c$-stratum $S^T_\beta$ is a linear subspace of $\PP_n$ and hence is 
nonsingular. Thus there is an obvious choice of lift $\widetilde{TG}^T_\beta$ to $X$ of
the Thom--Gysin map $TG^T_\beta$ which is given by the Thom--Gysin map associated
to the inclusion of $\bar{S}^T_\beta$ in $X$. We have
$$H^*_T(X) \cong H^*(BT)[\zeta]/\cali$$
where $\cali$ is the ideal in the polynomial ring $H^*(BT)[\zeta]$ generated
by the polynomial $(\zeta + \alpha_0) \cdots (\zeta + \alpha_n)$, while
$$H^*_T(\bar{S}^T_\beta) \cong H^*(BT)[\zeta] / \cali_\beta$$
where $\cali_\beta$ is the ideal generated by 
$$\prod_{\alpha_j.\beta \geq |\!|\beta|\!|^2} (\zeta + \alpha_j),$$
and
$$\widetilde{TG}^T_\beta ( \cali_\beta + \eta) = \cali + \eta \prod_{\alpha_j.\beta < |\!|\beta|\!|^2} (\zeta + \alpha_j)$$
where  
$$\prod_{\alpha_j.\beta < |\!|\beta|\!|^2} (\zeta + \alpha_j)$$
represents the equivariant Euler class $e_\beta^T$ of the normal bundle to $\bar{S}^T_\beta$ in $X$. Also
$$H^*_T(Z_\beta) \cong H^*(BT)[\zeta]/\calj_\beta$$
where $\calj_\beta$ is the ideal generated by
$$\prod_{\alpha_j.\beta = |\!|\beta|\!|^2} (\zeta + \alpha_j).$$
Hence by Lemma \ref{GG} there are lifts $\widetilde{TG}_\beta$ to $X$ of the Thom--Gysin maps
$$TG_\beta: H^{*-\lambda(\beta)}_{\stab(\beta)}(Z^{ss}_\beta) \cong H^{*-\lambda(\beta)}_G(S_\beta)
\to H^*_G(U_\beta)$$
represented by
$$\widetilde{TG}_\beta(\calj_\beta + \eta) = \cali + \frac{1}{|W|\cald} \sum_{w \in W} (-1)^w
w \left( \cald_\beta \; \eta \prod_{\alpha_j.\beta < |\!| \beta |\!|^2} (\zeta + \alpha_j) \right),$$
and by Lemma \ref{AA} the equivariant cohomology ring $H^*_G(X^{ss})$ is isomorphic
to the quotient of the polynomial ring $H^*(BT)[\zeta]$ by the ideal generated by
$(\zeta + \alpha_0) \cdots (\zeta + \alpha_n)$ 
and all polynomials in $\zeta$ with coefficients in $H^*(BT)$ of the form
$$\frac{1}{|W|\cald} \sum_{w \in W} (-1)^w
w \left( \cald_\beta \; \eta \prod_{\alpha_j.\beta < |\!| \beta |\!|^2} (\zeta + \alpha_j) \right)$$
for some $\beta \in \calb \setminus \{ 0 \}$ and $\eta \in H^*(BT)[\zeta]$.

When $G=SL(2;\CC)$ acts on $X = \PP_n$ as in Example \ref{example5.2} then the weights
are 
$$n\alpha, (n-2)\alpha, \ldots, -n \alpha$$
where $\alpha$ is a basis vector for $\liet$ and
$$\calb = \{ (2j-n) \alpha : j > \frac{n}{2} \} \cup \{ 0 \}.$$
Moreover $|W|=2$ and $\cald = 2 \alpha$, and if $\beta = (2j-n)\alpha \in \calb \setminus \{ 0 \}$
then $\cald_\beta = 1$ and
$\widetilde{TG}_\beta$ sends $ \calj_\beta + \eta(\zeta,\alpha)$ to
$$ \cali + \frac{1}{4\alpha}\left(
\eta(\zeta,\alpha) \prod_{k>j}(\zeta + (n-2k) \alpha)
- \eta(\zeta, - \alpha) \prod_{k>j} (\zeta - (n-2k) \alpha) \right).$$
Thus when $n$ is odd $H^*(X/\!/G)$ is generated as a $\QQ$-algebra by $\zeta$ and $\alpha$ with relations given
by
$$ \frac{1}{\alpha}\left(
\eta(\zeta,\alpha) \prod_{k>j}(\zeta + (n-2k) \alpha)
- \eta(\zeta, - \alpha) \prod_{k>j} (\zeta - (n-2k) \alpha) \right)$$
for all polynomials $\eta$ in $\zeta$ and $\alpha$.
\end{example}

\begin{example} \label{lastex}
When $G=SL(2;\CC)$ acts on $X=(\PP_1)^n$ the stratification $\{ \tilde{S}_{\ub}: \ub \in \tilde{\calb}  \}$
is described in Example \ref{previo}; it differs from the Morse stratification $\{ S_\beta: \beta \in \calb \}$
for $|\!|\mu |\!|^2$ only in that when $n$ is even the open stratum $S_0 = X^{ss}$ is decomposed into
the union of three strata, which are $\tilde{S}_{(0)} = X^s$ together with $\tilde{S}_{(T)}$ and $\tilde{S}_{(T,2)}$.
Here $H^*(BT) \cong \QQ[\alpha]$ where $\alpha$ has degree two and the nontrivial element of the Weyl 
group sends $\alpha$ to $-\alpha$, while $H^*_T(X)$ is generated by $n+1$ elements $\zeta_1, \ldots,
\zeta_n, \alpha$ of degree two subject to the relations
\begin{equation} \label{relations} (\zeta_1)^2 = \cdots = (\zeta_n)^2 = \alpha^2
\end{equation}
and $H^*_G(X)$ is generated by $\zeta_1, \ldots,
\zeta_n$ and $ \alpha^2$ subject to the same relations. The connected components $S_J$ of the strata
$S_\beta$ for $\beta \in \calb \setminus \{ 0 \}$ are indexed by subsets $J$ of  $\{1, \ldots, n\}$ of size
$|J| > n/2$, and their elements are sequences $(x_1, \ldots, x_n) \in (\PP_1)^n$ for which there is some
$p \in \PP_1$ satisfying $x_j = p$ if and only if $j \in J$. The connected components $S^T_J$ of the $T$-strata
$S^T_\beta$ are defined in the same way with $p=0$ or $p=\infty$. We have
$$\overline{S^T_J} \cong (\PP_1)^{n-|J|}$$
and if $\{1, \ldots, n\} \setminus J = \{ i_1, \ldots, i_{n-|J|} \}$ then the associated Thom--Gysin
map
$$\widetilde{TG}^T_J:H^{*-2|J|}_T(\overline{S^T_J}) \to H^*_T(X)$$
sends a polynomial
$p(\zeta_{i_1}, \ldots, \zeta_{i_{n-|J|}},\alpha) \in H^*_T((\PP_1)^{n-|J|})$ to
$$p(\zeta_{i_1}, \ldots, \zeta_{i_{n-|J|}},\alpha) \prod_{j \in J}(\zeta_j + \alpha) 
\in H^*_T((\PP_1)^{n}).$$
Thus by Lemma \ref{GG} a lift
$$\widetilde{TG}_J:H^{*-2|J|+2}_G({S_J}) \cong H^{*-2|J| + 2}(BT)  \to H^*_G(X)$$
of the Thom--Gysin map $TG_J$ is given by
$$\widetilde{TG}_J(p(\alpha)) = \frac{1}{4\alpha}\left( p(\alpha) \prod_{j \in J} (\zeta_j + \alpha)
- p(-\alpha) \prod_{j \in J} (\zeta_j - \alpha) \right).$$
It follows that $H^*_G(X^{ss})$ is generated by $\zeta_1, \ldots, \zeta_n$ and $\alpha^2$ subject to
the relations (\ref{relations}) and
\begin{equation} \label{rels}
\frac{1}{\alpha}\left(  \prod_{j \in J} (\zeta_j + \alpha)
-  \prod_{j \in J} (\zeta_j - \alpha) \right) = 0 =
 \prod_{j \in J} (\zeta_j + \alpha)
+ \prod_{j \in J} (\zeta_j - \alpha) 
\end{equation}
for all subsets $J$ of $\{ 1, \ldots , n \}$ with $|J|> n/2$.

When $n$ is even then the components of $\tilde{S}_{(T)}$ are indexed by partitions of
$\{ 1, \ldots, n\}$ into $J_1 \sqcup J_2$ where $|J_1| = |J_2| = n/2$, and  
their elements are sequences $(x_1, \ldots, x_n) \in (\PP_1)^n$ for which there are
$p_1 \neq p_2$ in $\PP_1$ satisfying $x_j = p_1$ if  $j \in J_1$ and $x_j = p_2$ if $j \in J_2$.
The components of $\tilde{S}_{(T,2)}$ are indexed by subsets $J_1$ of $\{1, \ldots, n\}$
with $|J_1|=n/2$, and 
 their elements are sequences $(x_1, \ldots, x_n) \in (\PP_1)^n$ for which there is some
$p_1 \in \PP_1$ satisfying $x_j = p_1$ if and only if $j \in J_1$, but no $p_2 \in \PP_1$
satisfying $x_j = p_2$ if $j \in \{ 1, \ldots, n \} \setminus J$. We have
$$H^*_G(S_{\{J_1,J_2\}}) \cong H^*(BT) \cong H^*_G(S_{J_1}\cup S_{\{J_1,J_2\}})$$
and the restriction map from $H^*_G(S_{J_1}\cup S_{\{J_1,J_2\}})$ to
$H^*_G(S_{J_1})$ is surjective.
Lifts to $H^*_G(X)$ of the Thom--Gysin maps
\begin{equation} \label{tgjj}
TG_{\{J_1,J_2\}}: H^{*-2n+4}_G(S_{\{ J_1, J_2\}}) \to H^*_G(X^{ss})
\end{equation}
and
\begin{equation} \label{tgj}
TG_{J_1}: H^{*-n+2}_G(S_{ J_1}) \to H^*_G(X^{ss}\setminus 
S_{\{ J_1, J_2\}})
\end{equation}
are given by
$$
TG_{\{J_1,J_2\}}\left(p(\alpha) \right) = $$
$$
\frac{1}{4\alpha}\left( p(\alpha) \alpha^{n/2 - 1} \prod_{j \in J_1} (\zeta_j + \alpha)
-  p(-\alpha)(-\alpha)^{n/2 -1} \prod_{j \in J_1} (\zeta_j - \alpha) \right)
$$
and
$$
TG_{J_1}\left (p(\alpha) \right) = \frac{1}{4\alpha} \left( p(\alpha)
 \prod_{j \in J_1} (\zeta_j + \alpha)
- p(-\alpha) \prod_{j \in J} (\zeta_j - \alpha)\right) .
$$
Thus the kernel of the restriction map
$$H^*_G(X) \to H^*_G(X^s)$$
is generated by the relations (\ref{rels}) for all subsets $J$ of $\{1, \ldots, n \}$ with $|J| \geq n/2$.
Note however that although the Thom--Gysin maps (\ref{tgjj}) and (\ref{tgj}) are injective,
the Thom--Gysin map associated to the inclusion of $S_{J_2}$ in
$X^{ss} \setminus(S_{\{J_1,J_2 \}} \cup S_{J_1})$ 
(which is just the composition of $TG_{J_2}$ as at (\ref{tgj}) with restriction from
$X^{ss} \setminus S_{\{J_1,J_2 \}} $ to
$X^{ss} \setminus(S_{\{J_1,J_2 \}} \cup S_{J_1})$ )
is not injective, and the restriction map from $H^*_G(X)$ to $H^*_G(X^s) \cong H^*(X^s/G)$ 
is not surjective when $n \geq 4$ is even. For example, when $n=4$ then
$X^s/G \cong \PP_1 \setminus \{ 0,1, \infty\}$ and so
$$\dim H^1(X^s/G) = 2$$
whereas the equivariant cohomology $H^*_G(X)$ of $X$ is all in even degrees.
\end{example}

\end{document}